%% file: AdaptiveFW.tex
\documentclass[twoside]{article}

\usepackage[accepted]{aistats2020}
\usepackage[round]{natbib}

\usepackage{empheq}
\usepackage[linesnumbered, ruled,vlined]{algorithm2e}
\usepackage{tikz}
\usepackage{multirow}
\usepackage{rotating}
\usepackage{xcolor, color, colortbl}
\usepackage{booktabs}
\usepackage{dsfont}
\usepackage{enumitem}
\usepackage{mdframed}
\usepackage{titlesec}
\usepackage{amsmath}
\usepackage{amsthm}
\usepackage{amssymb}
\usepackage{amsfonts}       
\usepackage{framed}
\usepackage{caption}
\usepackage{listings}
 
\definecolor{codegreen}{rgb}{0,0.6,0}
\definecolor{codegray}{rgb}{0.5,0.5,0.5}
\definecolor{codepurple}{rgb}{0.58,0,0.82}
\definecolor{backcolour}{rgb}{0.95,0.95,0.92}
 
\lstdefinestyle{mystyle}{
    backgroundcolor=\color{backcolour},   
    commentstyle=\color{codegreen},
    keywordstyle=\color{magenta},
    numberstyle=\tiny\color{codegray},
    stringstyle=\color{codepurple},
    basicstyle=\ttfamily\footnotesize,
    breakatwhitespace=false,         
    breaklines=true,                 
    captionpos=b,
    keepspaces=true,
    numbers=left,
    numbersep=5pt,
    showspaces=false,
    showstringspaces=false,
    showtabs=false,                  
    tabsize=2
}
 
\definecolor{shadecolor}{RGB}{242, 242, 242}

\definecolor{mydarkblue}{rgb}{0,0.08,0.45}
\usepackage[colorlinks=true,
    linkcolor=mydarkblue,
    citecolor=mydarkblue,
    filecolor=mydarkblue,
    urlcolor=mydarkblue,
    pdfview=FitH, draft=false]{hyperref}

\input defs.tex

\graphicspath{{./figures/}}

\newcommand{\tablefont}[1] {{\fontsize{8}{10}\sffamily{#1}}}

\begin{document}

%

%

\twocolumn[

  \aistatstitle{Linearly Convergent Frank-Wolfe with Backtracking Line-Search}

  \aistatsauthor{ Fabian Pedregosa \And Geoffrey N\'egiar \And  Armin Askari \And Martin Jaggi}

  \aistatsaddress{ Google Research \And University of California, Berkeley \And EPFL, Lausanne} ]

\begin{abstract}
Structured constraints in Machine Learning have recently brought the Frank-Wolfe (FW) family of algorithms back into the spotlight. While the classical FW algorithm has poor local convergence properties, Away-steps FW and Pairwise FW have emerged as improved variants with faster convergence. However, these improved variants suffer from two practical limitations: they require at each iteration to solve a 1-dimensional minimization problem to set the step-size and also require the Frank-Wolfe linear subproblems to be solved exactly. In this paper, we propose variants of Away-steps and Pairwise FW that lift both restrictions simultaneously. The proposed methods set the step-size based on a sufficient decrease condition, and do not require prior knowledge of the objective. Furthermore, they inherit all the favorable convergence properties of the exact line-search version, including linear convergence for strongly convex functions over polytopes. Benchmarks on different machine learning problems illustrate large performance gains of the proposed variants.
\end{abstract}

\section{Introduction}

The Frank-Wolfe (FW) or conditional gradient algorithm \citep{frank1956algorithm, levitin1966constrained, dem1967minimization} is a method for constrained optimization that solves problems of the form
\begin{empheq}[box=\mybluebox]{equation*}\tag{OPT}\label{eq:obj_fun}
  \minimize_{\xx \in \conv({\mathcal{A}})}\,\vphantom{\sum_i^n} f(\xx)\,,
\end{empheq}
where $f$ is a smooth function for which we have access to its gradient and $\conv(\mathcal{A})$ is the convex hull of $\mathcal{A}$; a bounded but potentially infinite set of elements in $\RR^p$ which we will refer to as \emph{atoms}.

The FW algorithm is one of the oldest methods for non-linear constrained optimization and has experienced a renewed interest in recent years due to its applications in machine learning and signal processing~\citep{jaggi2013revisiting}. Despite some favorable properties, the local convergence of the FW algorithm is known to be slow, achieving only a sublinear rate of convergence for strongly convex functions when the solution lies in the boundary~\citep{canon1968tight}.
To overcome these limitations, variants of the FW algorithms with better convergence properties have been proposed.
Two of these variants, the Away-steps FW~\citep{guelat1986some} and Pairwise FW~\citep{lacoste2015global} enjoy a linear rate of convergence over polytopes~\citep{lacoste2015global}.

Despite this theoretical breakthrough, Away-steps and Pairwise FW are not yet practical off-the-shelf solvers due to two main limitations.
The first and most important is that both variants rely on an \emph{exact line-search}. That is, they require to solve at each iteration \mbox{1-dimensional} subproblems of the form
\begin{equation}
  \label{eq:line-search}
  \argmin_{\gamma \in [0, \gamma_{\max}]} f(\xx_t + \gamma \dd_t)~,
\end{equation}
where $\dd_t$ is the update direction and $\gamma_{\max}$ is the maximum admissible step-size.
In a few cases like quadratic objectives, the exact line-search subproblem has a closed form solution. In most other cases, it is a costly optimization problem that needs to be solved at each iteration, making these methods impractical.
The second limitation is that they require access to an \emph{exact Linear Minimization Oracle} (LMO), which leaves out important cases like minimization over a trace norm ball where the LMO is computed up to some predefined tolerance.
In this paper we develop methods that lift both limitations simultaneously.


  { \begin{table*}[t]
      \centering
      \footnotesize
      \centering
      \footnotesize
      \setlength\tabcolsep{6pt}\begin{tabular}{c | c | c c c c c|}
        \cline{2-7}
        \multicolumn{1}{c|}{\multirow{6}{*}{\begin{sideways}\textbf{\sffamily Frank-Wolfe}\end{sideways}}}
                                                                         & \multirow{2}{*}{Related work}                                  & \multirow{1}{*}{non-convex}                                    & \multirow{1}{*}{approximate}                                   &
        \multirow{1}{*}{linear}                                          & \multirow{1}{*}{adaptive}                                      & \multirow{1}{*}{bounded}                                                                                                                                                                                                       \\
        \multicolumn{1}{c|}{}                                            &                                                                & \multirow{-1}{*}{analysis}                                     & \multirow{-1}{*}{subproblems}                                  & \multirow{-1}{*}{convergence} & \multirow{-1}{*}{step-size} & \multirow{-1}{*}{backtracking} \\
        \cline{2-7}

                                                                         & \multicolumn{1}{c|}{\emph{This work}{\cellcolor{Gray}}}        &
        \multirow{-1}{*}{\color{mygreen}\large\cmark} \cellcolor{Gray}
        \cellcolor{Gray}                                                 &
        \cellcolor{Gray} \multirow{-1}{*}{\color{mygreen}\large\cmark}   & \cellcolor{Gray} \multirow{-1}{*}{\color{mygreen}\large\cmark} & \cellcolor{Gray} \multirow{-1}{*}{\color{mygreen}\large\cmark} & \cellcolor{Gray} \multirow{-1}{*}{\color{mygreen}\large\cmark}                                                                                                \\

                                                                         &
        \multirow{-1}{*}{\citep{lacoste2015global}}                      & \multirow{-1}{*}{\red\large\xmark}                             & \multirow{-1}{*}{\red\large\xmark}                             &
        \multirow{-1}{*}{\color{mygreen}\large\cmark}
                                                                         &
        \multirow{-1}{*}{\red\large\xmark}                               &
        N/A                                                                                                                                                                                                                                                                                                                                                                \\

        \multicolumn{1}{c|}{\multirow{4}{*}}                             &
        \multirow{-1}{*}{\citep{beck2015cyclic}}\cellcolor{Gray}         &
        \cellcolor{Gray} \red\large\xmark                                & \cellcolor{Gray} {\color{mygreen}{\large\cmark}}$^\dag$        &
        \multirow{-1}{*}{\red\large\xmark} \cellcolor{Gray}
        \cellcolor{Gray}                                                 &
        \cellcolor{Gray} \multirow{-1}{*}{\color{mygreen}{\large\cmark}} & \cellcolor{Gray} \multirow{-1}{*}{\red\large\xmark}                                                                                                                                                                                                                                             \\

        \multicolumn{1}{c|}{\multirow{4}{*}}                             &
        \multirow{-1}{*}{\citep{dunn1980convergence}}                    &
        \multirow{-1}{*}{\color{mygreen}{\large\cmark}}                  & \multirow{-1}{*}{\red\large\xmark}
                                                                         & \multirow{-1}{*}{\red\large\xmark}
                                                                         & \multirow{-1}{*}{\color{mygreen}{\large\cmark}}                & \multirow{-1}{*}{\red\large\xmark}                                                                                                                                                                                             \\
        \cline{2-7}\addlinespace[0.2cm]

        \cmidrule{2-7}\vspace{-3.2ex}                                                                                                                                                                                                                                                                                                                                      \\
        \multicolumn{1}{c|}{\multirow{2}{*}{\begin{sideways}\textbf{\sffamily MP}\end{sideways}}}
                                                                         & \multirow{-1}{*}{\emph{This work}}
        \cellcolor{Gray}                                                 &
        \multirow{-1}{*}{\color{mygreen}\large\cmark} \cellcolor{Gray}
        \cellcolor{Gray}                                                 &
        \cellcolor{Gray} \multirow{-1}{*}{\color{mygreen}\large\cmark}   & \cellcolor{Gray} \multirow{-1}{*}{\color{mygreen}\large\cmark} & \cellcolor{Gray} \multirow{-1}{*}{\color{mygreen}\large\cmark} & \cellcolor{Gray} \multirow{-1}{*}{\color{mygreen}\large\cmark}                                                                                                \\

        \multicolumn{1}{c|}{\multirow{4}{*}}                             &
        \multirow{-1}{*}{\citep{locatello17a}}                           & \multirow{-1}{*}{\red\large\xmark}                             & \multirow{-1}{*}{\color{mygreen}\large\cmark}                  &
        \multirow{-1}{*}{\color{mygreen}\large\cmark}
                                                                         &
        \multirow{-1}{*}{\red\large\xmark}                               &
        N/A                                                                                                                                                                                                                                                                                                                                                                \\ \cline{2-7}
      \end{tabular}
      \caption{{\bfseries Comparison with related work.} \emph{non-convex analysis}: convergence guarantees for problems with a non-convex objective.
        \emph{approximate subproblems}: convergence guarantees cover the case in which linear subproblems are solved approximately.
        \emph{linear convergence}: guaranteed linear rate of convergence (under hypothesis).
        \emph{adaptive step-size}: step-size is set using local information of the objective.
        \emph{bounded backtracking}: explicit bound for the total number of inner iterations in adaptive step-size methods. $^\dagger=$ assumes cartesian product domain. 
      }\label{table:related_work}
    \end{table*}
  }


Our \textbf{main contribution} is the design and analysis of variants of Away-steps and Pairwise FW that $i)$ don't require access to an exact line-search or knowledge of properties of the objective like its curvature or Lipschitz constant, and $ii)$ admits the FW subproblems to be solved approximately. We describe our approach in \S \ref{section:methods}. Although our main motivation is to develop practical variants of Away-steps and Pairwise FW, we also show that this technique extends to other methods like FW and Matching Pursuit.

We develop in \S\ref{section:theory} a convergence rate analysis for the proposed methods. The obtained rates match asymptotically the best known bounds on convex, strongly convex and non-convex problems, including linear convergence for strongly convex functions.

Finally, we show in \S\ref{section:experiments} benchmarks between the proposed and related methods, and discuss the importance of large step-sizes in Pairwise FW.

\subsection{Related work}

We comment on the most closely related ideas, summarized in Table~\ref{table:related_work}.

Away-Steps FW~\citep{guelat1986some} is a popular variant of FW that adds the option to move away from an atom in the current representation of the iterate. In the case of a polytope domain, it was recently shown to enjoy a linear convergence rate for strongly convex objectives \citep{garber2013linearly,lacoste2015global}.
Pairwise FW \citep{lacoste2015global} simplifies the above-described variant by replacing the two kinds of steps by a single step modifying the weights of only two atoms. It generalizes the algorithm of \citet{mitchell1974finding} used in geometry and SMO \citep{platt1998sequential} for training SVMs. These methods all require exact line-search.

Variants of FW that, like the proposed methods, set the step-size based on a local decrease condition have been described by \citet{dunn1980convergence} and \citet{beck2015cyclic}, but none of these methods achieve a linear convergence rate to the best of our knowledge.

Matching Pursuit (MP)~\citep{mallat1993matching} is an algorithm for constrained optimization problems of the form \eqref{eq:obj_fun} with $\conv(\AA)$ replaced by $\text{linspan}(\AA)$, the linear span of $\AA$.
\citet{locatello2018revisiting} has recently shown that MP and FW are deeply related. 
We show that our algorithm and convergence results also extend naturally to MP, and as a byproduct of our analysis we obtain the first convergence rate for MP on non-convex objectives to the best of our knowledge.


\paragraph{Notation.} Throughout the paper we denote vectors and vector-valued functions in lowercase boldface (e.g. $\xx$ or $\argmin$), matrices in uppercase boldface letters (e.g. $\boldsymbol D$), and sets in caligraphic letters (e.g., $\mathcal{A}$).
We say a function $f$ is $L$-smooth if it is differentiable and its gradient is $L$-Lipschitz continuous, that is, if it verifies $\|\nabla f(\xx) - \nabla f(\yy)\| \leq L\|\xx - \yy\|$ for all $\xx, \yy$ in the domain.
A function is { $\mu$-strongly convex} if $f - \frac{\mu}{2}\|\cdot\|^2$ is convex. $\|\cdot\|$ denotes the euclidean norm.

\section{Methods}
\label{section:methods}

In this section we describe the core part of our contribution, which is a strategy to select the step-size in FW-type algorithms.

Since this strategy can be applied very broadly to Frank-Wolfe variants including Away-steps, Pairwise, classical FW and Matching Pursuit, we describe it within the context of a generic FW-like algorithm.
This generic algorithm is detailed in Algorithm~\ref{alg:metafw} and depends on two key functions: \texttt{update\_direction} and \texttt{step\_size}. The first one computes the direction that we will follow to compute the next iterate and its implementation will depend on the FW variant. The second one will choose an appropriate step-size based upon local information of the objective and is the key novelty of this algorithm.
We now describe them in more detail.

\begin{algorithm}[t]

  {\bfseries Input:} $\xx_0 \in \text{conv}(\mathcal{A})$, initial Lipschitz estimate $L_{-1} > 0$, tolerance $\varepsilon \geq 0$, subproblem quality $\delta \in (0, 1]$


  \For{$t=0, 1 \ldots $}{
      $\dd_t, \gamma_t^\text{max} = \text{update\_direction}(\xx_t, \nabla f_t)$\label{l:get_direction}

    $g_t = \langle -\nabla f(\xx_t), \dd_t \rangle$

    \lIf{ ${g}_t \leq \delta\varepsilon\,$}{ \Return $\xx_t$}

    \begin{shaded}$\gamma_t, L_{t} = \text{step\_size}(f, \dd_t, \xx_t, g_t, L_{t-1}, \gamma_t^\text{max})$\end{shaded}
      
    $\xx_{t+1} = \xx_t + \gamma_t \dd_t$

  }
  \caption{Generic FW with backtracking}\label{alg:metafw}
\end{algorithm}


\subsection{Update direction}

In this subsection we describe \texttt{update\_direction} in Algorithm~\ref{alg:metafw}, the function that computes the update direction $\dd_t$ and the maximum allowable step-size~$\gamma_t^{\max}$.
While its implementation varies according to the FW variant, all of them require to solve one or two linear problems, often referred to as linear minimization oracle (LMO).

The first of these subproblems is the same for all variants and consists in finding $\ss_t$ in the domain such that:
\begin{equation}\label{eq:lmo}
  \langle \nabla f(\xx_t),\ss_t - \xx_t\rangle \leq \delta\min_{\ss \in \mathcal{A}}\langle \nabla f(\xx_t),  \ss - \xx_t \rangle~.
\end{equation}
Here, we introduce a \emph{subproblem quality} parameter $\delta \in (0, 1]$  that allows this subproblem to be solved approximately. When $\delta = 1$, the problem is solved exactly and becomes $\argmin_{\ss \in \mathcal{A}} \langle \nabla f(\xx_t), \ss \rangle$, which consists in selecting the atom that correlates the most with the steepest descent direction, $-\nabla f(\xx_t)$.

Away-steps and Pairwise FW will also require to solve another linear subproblem, this time over the \emph{active set} $\mathcal{S}_t$. This is the set of atoms with non-zero weight in the decomposition of $\xx_t$.
More formally, the active set $\mathcal{S}_t \subseteq \mathcal{A}$ is the set of atoms that have non-zero weight $\balpha_{\ss, t} > 0$ in the expansion $\xx_t = \sum_{\ss \in \mathcal{S}_t} \balpha_{\vv, t} \ss$.

The linear subproblem that needs to be solved consists in finding $\vv_t$ such that:
\begin{equation}\label{eq:away_lmo}
  \langle \nabla f(\xx_t),  \xx_t - \vv_t \rangle \leq \delta\min_{\vv \in {S}_t}\langle \nabla f(\xx_t), \xx_t - \vv\rangle~.
\end{equation}
Unlike the previous linear subproblem, this time the problem is over the typically much smaller active set~$\mathcal{S}_t$.
As before, $\delta \in (0, 1]$ allows this subproblem to be solved approximately. 
When $\delta=1$, the subproblem becomes $\argmax_{\vv \in \mathcal{S}_t} \langle \nabla f(\xx_t), \vv\rangle$, which can be interpreted as selecting the atom in the active set that correlates the most with the steepest ascent direction~$\nabla f(\xx_t)$.


FW, AFW and PFW then combine the solution to these linear subproblems in different ways, and Line~\ref{l:get_direction}'s \texttt{update\_direction} is implemented as:

\begin{itemize}[leftmargin=*]
  \item FW returns $\dd_t = \ss_t - \xx_t$ and $\gamma_t^{\max}=1$: the next iterate will be a convex combination of $\xx_t$ and $\ss_t$.
  \item AFW considers directions $\ss_t - \xx_t$ and $\xx_t - \vv_t$, and chooses the one that correlates the most with $-\nabla f(\xx_t)$. $\gamma_t^{\max} = 1$ if $\dd_t = \ss_t - \xx_t$ and $\alpha_{\vv_t} / ( 1- \alpha_{\vv_t})$ otherwise, where $\alpha_{\vv_t}$ is the weight associated with~$\vv_t$ in the decomposition of $\xx_t$ as a convex combination of atoms.
  \item PFW uses $\dd_t = \ss_t - \vv_t$, shifting weight from $\vv_t$ to~$\ss_t$ in our current iterate, and $\gamma_t^{\max} = \alpha_{\vv_t}$.
  \item MP uses $\dd_t = \ss_t$ and $\gamma^\text{max} = +\infty$, since the constraint set is not bounded.

\end{itemize}


\subsection{Backtracking line-search}

In this subsection we describe the step-size selection routine \texttt{step\_size} (Algorithm~\ref{alg:linesearch}). 
This is the main novelty in the proposed algorithms, and allows for the step-size to be computed using only \emph{local} properties of the objective, as opposed to other approaches that use global quantities like the gradient's Lipschitz constant. As we will see in \S\ref{section:experiments}, this results in step-sizes that are often more than an order of magnitude larger than those estimated using global quantities.  

\begin{algorithm}[t]
  \SetAlgoLined\DontPrintSemicolon
  \SetKwFunction{algo}{algo}\SetKwFunction{proc}{step\_size}
  \SetKwProg{myalg}{Algorithm}{}{}
  \SetKwProg{myproc}{Procedure}{}{}
  \myproc{\proc{$f$,$\dd_t$, $\xx_t$, $g_t$, $L_{t-1}$,  $\gamma_{\max}$}}{

  Choose $\tau > 1$,  $\eta \leq 1$
  
  Choose $M \in [\eta L_{t-1}, L_{t-1}]$






  $\gamma = \min\left\{{{g}_t}/{(M\|\dd_t\|^{2})}, \gamma_{\max}\right\}$\label{l:step_size}


  \While{$f(\xx_t + \gamma \dd_t) > Q_t(\gamma, M) $ \label{l:afw_while}}{
  $M = \tau M$ \label{l:increase_tau}

  $\gamma = \min\left\{{{g}_t}/{(M\|\dd_t\|^{2})}, \gamma_{\max}\right\}$\label{l:inner_linesearch}


  }

  \KwRet $\gamma$, $M$\;}
  \caption{Backtracking for FW variants}\label{alg:linesearch}
\end{algorithm}

Minimizing the exact line-search objective $\gamma \mapsto {f(\xx_t + \gamma \dd_t)}$ yields the highest decrease in objective but can be a costly optimization problem. To overcome this, we will replace the exact line-search objective with the following quadratic approximation:
\begin{equation}\label{eq:surrogate}
  Q_t(\gamma, M) = f(\xx_t) - \gamma g_t + \frac{\gamma^2 M}{2}\|\dd_t\|^2~.
\end{equation}
This approximation has the advantage that its minimum over $\gamma \in [0, \gamma^{\text{max}}]$ can be computed in closed form, which gives the step-size used in line \ref{l:step_size}:
\begin{equation}
  \gamma^\star_M = \min\left\{\frac{{g}_t}{M\|\dd_t\|^{2}}, \gamma^\text{max}\right\}~,
\end{equation}

The quality of this quadratic approximation will depend on the \emph{Lipschitz estimate} parameter $M$. This parameter needs to be carefully selected to maintain the convergence guaratees of exact line-search, while keeping the number of objective function evaluations to a minimum.

This is achieved through the strategy implemented in Algorithm~\ref{alg:linesearch}. The algorithm initializes the Lipschitz estimate $M$ to a value between $\eta L_{t-1}$ and the previous iterate $L_{t-1}$, where $\eta$ is a user-defined parameter (default values discussed later). A value of $\eta=1$ is admissible but would not allow the Lipschitz estimate to decrease through the optimization, and we have observed empirically a drastic benefit in doing so.

The algorithm then defines a candidate step-size $\gamma$ (Line~\ref{l:step_size}) and checks whether the following \emph{sufficient decrease} condition is verified for this step-size
\begin{equation*}
  f(\xx_t + \gamma \dd_t) \leq Q_t(\gamma, M)~, \gamma\!=\!\min\left\{{{g}_t}/{(M\|\dd_t\|^{2})}, \gamma^\text{max}\right\}.
\end{equation*}
If it is not verified, we increase this constant by a power factor of $\tau > 1$ (Line \ref{l:increase_tau}). By the properties of $L$-smooth functions, we know that this condition is verified for all $M \geq L$, and so this loop has a finite termination.

Once this condition is verified, the current step-size is accepted and the value of $M$ is assigned the name $L_t$.
\begin{minipage}[ht]{.43\linewidth}
  \vspace{0.0em}
  Geometrically, the sufficient decrease condition ensures that the quadratic approximation is an upper bound at its constrained minimum of the line-search objective. We emphasize that this upper
  \unskip\parfillskip 0pt \par
\end{minipage}
\begin{minipage}[ht]{.51\linewidth}
  \hspace*{0.5em}{\begin{tikzpicture}[domain=-0.1:3, baseline=(current bounding box.north)]
      \draw[->] (-0.2,0) -- (2.8,0) node[below] {};
      \draw[->] (0,-0.2) -- (0,2.9) node[above] {};
      \node[circle,fill=black,inner sep=0pt,minimum size=3pt,label=below:{$\gamma\!=\!0$}] (a) at (0.,0) {};
      \node[circle,fill=black,inner sep=0pt,minimum size=3pt,label=below:{$\gamma\!=\!\gamma_t^{\max}$}] (a) at (2.4,0.) {};
      \node[circle,fill=black,inner sep=0pt,minimum size=3pt,label=below:{$\color{brown}\gamma_{t}$}] (a) at (1.17893,0) {};
      \node[circle,draw=black, fill=white, inner sep=0pt,minimum size=5pt] (b) at (0. ,2.4) {};
      \node[circle,draw=black, fill=white, inner sep=0pt,minimum size=5pt] (b) at (2.4 ,0.52) {};
      \draw[gray, dashed] (2.4,0)--(2.4,0.45);
      \draw[brown, dashed] (1.17893,0)--(1.17893,1.3675);
      \draw[color=brown,densely dashed, domain=0.0:2.4, line width=0.4mm] plot ({\x}, {2.41 - 1.7684 * (\x) + 0.75* (\x) * (\x)}) node[right, align=left] {$Q_{t}(\gamma, L_t)$};;
      \draw[color=black, domain=0.:2.4, line width=0.5mm] plot (\x, { 2.7 * exp(-0.8*\x-0.2) +0.2})
      node[right] {$f(\xx_t\!\!+\!\!\gamma \dd_t)$};
    \end{tikzpicture}}
  \label{fig:sufficient_decrease}
\end{minipage}

\vspace{-0.25em}
bound does not need to be a global one,  as it only holds at $\gamma_t$.
This allows for smaller $L_t$ than the global Lipschitz constant $L$, and therefore larger step-sizes. As we will see in \S\ref{section:theory}, this translates into faster convergence rates that depend on $L_t$, as well as faster empirical convergence (\S\ref{section:experiments}).

\paragraph{Default and initial parameters.}
Algorithm \ref{alg:metafw} requires an (arbitrary) initial value for the Lipschitz estimate $L_{-1}$.
We found the following heuristic using the definition of Lipschitz continuity of the gradient to work well in practice. Select a small constant $\varepsilon$, say $10^{-3}$, and compute
$L_{-1} = {\|\nabla f(\xx_0) - \nabla f(\xx_0 + \varepsilon \dd_0)\|}/{(\varepsilon\|\dd_0\|)}$.

The \texttt{step\_size} depends on hyperparameters $\eta$ and $\tau$. Although the algorithm is guaranteed to converge for any  $\eta \leq 1$, $\tau > 1$, we recommend $\eta = 0.9, \tau=2$, as we found that it performs well in a variety of scenarios. These are the values used throughout benchmarks \S\ref{section:experiments}.

This method also requires to choose the initial value of the Lipschitz estimate $M$ to a value between $\eta L_{t-1}$ and $L_{t-1}$. A choice that we found to work remarkably well in practice is to initialize it to 
\begin{equation}
  M = \text{clip}_{[\eta L_{t-1}, L_{t-1}]}\left(\frac{g_t^2}{2 (f_{t-1} - f_t) \|\dd_t\|^2}\right)~.
\end{equation}
The value inside the clip function corresponds to the optimal value of $M$ for a quadratic interpolation between the previous two iterates and the derivative of the line-search objective $f(\xx_t + \gamma \dd_t)$ at $\gamma=0$. Since this value might be outside of the interval $[\eta L_{t-1}, L_{t-1}]$, we clip the result to this interval. 

\paragraph{Pseudocode and implementation details.}
A practical implementation of these algorithms depends on other details that are not specific to the backtracking variant, such as efficiently maintaining the active-set in the case of Away-steps and Pairwise. For completeness, \ref{apx:pseudocode} contains a full pseudocode for all these algorithms. A Python implementation of these methods, as well as the benchmarks used in \S\ref{section:experiments} will be made open source upon publication of this manuscript. 

\section{Analysis}
\label{section:theory}

In this section, we provide a convergence rate analysis of the proposed methods, showing that all  enjoy a $\mathcal{O}(1/\sqrt{t})$ convergence rate for non-convex objectives (Theorem~\ref{thm:noncvx}), a stronger $\mathcal{O}(1/t)$ convergence rate for convex objectives (Theorem~\ref{thm:convex}), and for some variants linear convergence for strongly convex objectives over polytopes (Theorem~\ref{thm:linear_conv_afw}).

\vspace{-0.5em}\paragraph{Notation.} In this section we make use of the following extra notation:
\vspace{-0.5em}\begin{itemize}[leftmargin=*]
  \item For convenience we will refer to the variants of FW, Away-steps FW, Pairwise FW and MP with backtracking line-search as \hyperref[alg:afw]{AdaFW}, \hyperref[alg:aafw]{AdaAFW}, \hyperref[alg:apfw]{AdaPFW} and \hyperref[alg:amp]{AdaMP} respectively.
  \item  We denote the \emph{objective suboptimality} at step $t$ as $h_t = f(\xx_t) - \min_{\xx \in \conv(\mathcal{A})}f(\xx)$.
  \item \emph{Good and bad steps}. Our analysis, as that of \citet{lacoste2015global}, relies on a notion of ``good'' and ``bad'' steps. We define bad steps as those that verify $\gamma_t = \gamma_t^{\max}$ and $\gamma_t^{\max} < 1$ and good steps as any step that is not a bad step. The name ``bad steps'' makes reference to the fact that we won't be able to bound non-trivially the improvement for these steps. For these steps we will only be able to guarantee that the objective is non-increasing. AdaAFW and AdaPFW both may have bad steps. Let us denote by $N_t$ the number of ``good steps'' up to iteration $t$. We can lower bound the number of good steps by
        \begin{align}\label{eq:bad_steps}
          N_t \geq t/2 \text{ for  \hyperref[alg:aafw]{AdaAFW}}\,, \\
          N_t \geq t/(3|\mathcal{A}|! + 1) \text{ for \hyperref[alg:apfw]{AdaPFW}}
        \end{align}
        where it is worth noting that the last bound for \hyperref[alg:apfw]{AdaPFW} requires the set of atoms $\mathcal{A}$ to be finite. The proof of these bounds can be found in \ref{scs:bound_bad_steps} and are a direct translation of those in \citep{lacoste2015global}. We have found these bounds to be very loose, as in practice the fraction of bad/good steps is negligible, commonly of the order of $10^{-5}$ (see last column of the table in Figure \ref{fig:main_figure}).

  \item \emph{Average and maximum of Lipschitz estimates}. In order to highlight the better convergence rates that can be obtained by adaptive methods we introduce the average and maximum estimate over good step-sizes. Let $\mathcal{G}_t$ denote the indices of good steps up to iteration $t$. Then we define the average and maximum Lipschitz estimate as
        \begin{align}
          \overline{L}_t   & \defas \frac{1}{N_t}{\textstyle\sum_{k \in \mathcal{G}_t}} L_k \\
          \quad L_t^{\max} & \defas \textstyle{\max_{k \in \mathcal{G}_t}} L_k
        \end{align}
        respectively. 
        In the worst case, both quantities can be upper bounded by $\max\{\tau L, L_{-1}\}$  (Proposition~\ref{prop:lipschitz_estimate_bound}), which can be used to obtain asymptotic convergence rates. This bound is however very pessimistic. We have found that in practice $\overline{L}_t$ is often more than 100 times smaller than $L$ (see second to last column of the table in Figure \ref{fig:main_figure}).
\end{itemize}

{ \begin{table}[t]
  \centering
  \setlength\tabcolsep{4pt}\begin{tabular}{c | c | c c c|}
    \cline{2-5}
    \multicolumn{1}{c|}{}
                          & \multirow{2}{*}{\sffamily Algorithm}                                   & \multirow{2}{*}{\sffamily Non-convex\!\!} & \multirow{2}{*}{\sffamily Convex} &
    \multirow{1}{*}{\sffamily Strongly}                                                                                                                                              \\
    \multicolumn{1}{c|}{} &                                                              &                             &                         & \multirow{1}{*}{\sffamily convex} \\
    \cline{2-5}

                          &                                                              &                             &                         &                         \\
                          & \multirow{-2}{*}{AdaAFW}                                     &
    \multirow{-2}{*}{\large$\mathcal{O}(\frac{1}{\delta\sqrt{t}})$}
                          & \multirow{-2}{*}{\large$\mathcal{O}(\frac{1}{\delta^2{t}})$}
                          & \multirow{-2}{*}{\large$\mathcal{O}({(1\!-\!\delta^2\rho)}^t)$}                                                                                                  \\

                          &                                                              &                             &                         &                         \\
                          & \multirow{-2}{*}{AdaPFW}                                     &
    \multirow{-2}{*}{\large$\mathcal{O}(\frac{1}{\delta\sqrt{t}})$}
                          & \multirow{-2}{*}{\large$\mathcal{O}(\frac{1}{\delta^2{t}})$}
                          & \multirow{-2}{*}{\large$\mathcal{O}((1\!-\!\delta^2\rho)^t)$}                                                                                                  \\

                          &                                                              &                             &                         &                         \\
                          & \multirow{-2}{*}{AdaFW}                                      &
    \multirow{-2}{*}{\large$\mathcal{O}(\frac{1}{\delta\sqrt{t}})$}
                          & \multirow{-2}{*}{\large$\mathcal{O}(\frac{1}{\delta^2{t}})$}
                          & \multirow{-2}{*}{\large$\mathcal{O}(\frac{1}{\delta^2{t}})$}                                                                                   \\

                          &                                                              &                             &                         &                         \\
                          & \multirow{-2}{*}{AdaMP}                                      &
    \multirow{-2}{*}{\large$\mathcal{O}(\frac{1}{\delta\sqrt{t}})$}
                          & \multirow{-2}{*}{\large$\mathcal{O}(\frac{1}{\delta^2{t}})$}
                          & \multirow{-2}{*}{\large$\!\mathcal{O}((1\!-\!\delta^2\rho_{\text{MP}})^{t})\!$}                                                                                                  \\

    \cline{2-5}
  \end{tabular}
  \caption{{\bfseries Convergence rate summary} on non-convex, convex and strongly convex objectives. For non-convex objectives, bound is on the minimum FW gap (MP gap in the case of AdaMP), in other cases its on the objective suboptimality.
  }\label{table:related_work2}
\end{table}
}

Our new convergence rates are presented in the following theorems, which consider the cases of  non-convex, convex and strongly convex objectives. The results are discussed in \S\ref{scs:discussion}  and the proofs can be found in \ref{apx:proof_nonconvex}, \ref{apx:proof_convex} and \ref{apx:linear_proofs} respectively.

\subsection{Overhead of backtracking}

Evaluation of the sufficient decrease condition Algorithm~\ref{alg:linesearch} requires two extra evaluations of the objective function.
If the condition is verified, then it is only evaluated at the current and next iterate. FW requires anyway to compute the gradient at these iterates, hence in cases in which the objective function is available as a byproduct of the gradient this overhead becomes negligible.

Furthermore, we can provide a bound on the total number of evaluations of the sufficient decrease condition:
\begin{theorem}\label{thm:bound_backtracking}
  Let $n_t$ be the total number of evaluations of the sufficient decrease condition up to iteration $t$. Then we have
  \begin{equation*}
    n_t \leq \left[ 1 -\frac{\log \eta}{\log \tau}\right](t+1) + \frac{1}{\log\tau}\max\left\{\log\frac{\tau L}{L_{-1}}, 0\right\}~,
  \end{equation*}
\end{theorem}

This result highlights the trade-off faced when choosing $\eta$.
Minimizing it with respect to $\eta$ gives $\eta=1$, in which case $(1 - \log\eta/\log\tau)=1$ and so there's an asymptotically vanishing number of failures in the sufficient decrease condition. Unfortunately, $\eta=1$ also forbids the Lipschitz estimate to decrease along the optimization. Ideally, we would like $\eta$ small enough so that the Lipschitz estimate decreases when it can, but not too small so that we waste too much time in failed sufficient decrease evaluations.

As mentioned before, we recommend parameters $\eta = 0.9$, $\tau=2$. With these values, we have that $\left[ 1 -\frac{\log \eta}{\log \tau}\right] \leq 1.16$, and so asymptotically no more than 16\% of the iterates will result in more than one evaluations of the sufficient decrease condition.

\subsection{Non-convex objectives}

\paragraph{Gap function.} Convergence rates for convex and strongly convex functions are given in terms of the objective function suboptimality or a primal-dual gap. As the gap upper-bounds (i.e. certifies) the suboptimality, the latter is a stronger result in this scenario.
In the case of non-convex objectives, as is common for first order methods, we will only be able to guarantee convergence to a stationary point, defined as any element $\xx^\star \in \conv(\mathcal{A})$ such that $\langle \nabla f(\xx^\star), \xx - \xx^\star \rangle \geq 0~\text{ for all $\xx \in \conv(\mathcal{A})$}$~\citep{bertsekas1999nonlinear}.

Following \citet{lacoste2016convergence, reddi2016stochastic}, for FW variants we will use as convergence criterion the FW gap, defined as $g^{\text{FW}}(\xx) = {\max_{\ss \in \conv(\mathcal{A})} \langle \nabla f(\xx), \xx - \ss\rangle}$. From the definition of stationary point it is clear that the FW gap is zero only at a stationary point. These rates are also valid for AdaMP, albeit for the more appropriate gap function $g^{\text{MP}}$ detailed in \ref{apx:proof_nonconvex}.

\begin{theorem}\label{thm:noncvx}
  Let $\xx_{t}$ denote the iterate generated by any of the proposed algorithms after $t$ iterations, with $N_{t+1} \geq 1$.
  Then  we have:
  \begin{align}\label{eq:rate_nonconvex}
    &~ \lim_{t\to \infty}g(\xx_t) = 0\qquad \text{ and } \\ &\min_{k=0,\dots, t} g(\xx_k) \leq \frac{C_t}{\delta\sqrt{N_{t+1}}} = \mathcal{O}\left(\frac{1}{\delta\sqrt{t}}\right),
  \end{align}
  where $C_t = \max\{2 h_0, L_t^{\max}\diam(\AA)^2\}$ and $g = g^{\text{FW}}$ is the FW gap for \hyperref[alg:afw]{AdaFW}, \hyperref[alg:aafw]{AdaAFW}, \hyperref[alg:apfw]{AdaPFW} and $C_t =  \radius(\mathcal{A})\sqrt{2 h_0{\overline{L}_{t+1}}}$ and $g = g^{\text{MP}}$ is the MP gap for \hyperref[alg:amp]{AdaMP}.

\end{theorem}

\subsection{Convex objectives}

For convex objectives we will be able to improve the results of Theorem~\ref{thm:noncvx}. We will first state the convergence results for FW variants and then for MP.

For adaptive FW variants, we will be able to give an $\mathcal{O}(1/\delta^2 t)$ convergence rate on the primal-dual gap, which trivially implies a bound on the objective suboptimality. In order to define the primal-dual gap, we define the following \emph{dual} objective function
\begin{equation}
  \psi(\uu) \defas -f^*(\uu) - \sigma_{\conv(\mathcal{A})}(-\uu)~,
\end{equation}
where $f^*$ denotes the convex conjugate of $f$ and $\sigma_{\conv(\mathcal{A})}(\xx) \defas \sup\{ \langle \xx, \aa \rangle: \aa\in \conv(\mathcal{A})\}$ is the support function over $\conv(\mathcal{A})$, which is the convex conjugate of the indicator function. Note that $\psi$ is concave and that when $f$ convex, we have by duality $\min_{\xx \in \conv(\mathcal{A})}f(\xx_t) = \max_{\uu \in \RR^p} \psi(\uu)$.

\begin{theorem}[FW variants]\label{thm:convex}
  Let $f$ be convex, $\xx_{t}$ denote the iterate generated by any of the proposed FW variants (\hyperref[alg:afw]{AdaFW}, \hyperref[alg:aafw]{AdaAFW}, \hyperref[alg:apfw]{AdaPFW}) after $t$ iterations, with $N_t \geq 1$, and let $\uu_t$ be defined recursively as $\uu_0 = \nabla f(\xx_0)$, $\uu_{t+1} = (1 - \xi_t)\uu_{t} + \xi_t\nabla f(\xx_t)$, where $\xi_t = 2 / (\delta N_t + 2)$ if $t$ is a good step and $\xi_t = 0$ otherwise. Then we have:
  \begin{align}
    h_t & \leq f(\xx_t) - \psi(\uu_t)                                                                                                                                            \\
        & \quad\leq \frac{2\overline{L}_{t}\diam(\AA)^2}{ \delta^2 N_t + \delta} + \frac{2(1 - \delta )}{\delta^2 N_t^2 + \delta N_t} \big(f(\xx_0) - \psi(\uu_0) \big)\nonumber \\
        & \quad= \mathcal{O}\left(\frac{1}{\delta^2 t}\right).
  \end{align}
\end{theorem}

\subsection{Strongly convex objectives}

\newcommand{\TheoremLinear}{
  Let $f$ be $\mu$--strongly convex. Then for \hyperref[alg:aafw]{AdaAFW}, \hyperref[alg:apfw]{AdaPFW} or \hyperref[alg:amp]{AdaMP} we have the following linear decrease for each good step $t$:
  \begin{equation}
    h_{t+1} \leq (1 - \delta^2 \rho_t) h_t,\vspace{-2mm}
  \end{equation}
  where
  \begin{flalign*}
    \rho_t &=  \dfrac{\mu}{4{L}_t}\left(\dfrac{\PWidth(\AA)}{\diam(\AA)}\right)^2    \text{for \hyperref[alg:aafw]{AdaAFW} and \hyperref[alg:apfw]{AdaPFW}},& \\
    \rho_t &=  \dfrac{\mu}{{L}_t} \left(\dfrac{\mDW(\AA)}{\radius(\AA)}\right)^2 \text{ for \hyperref[alg:amp]{AdaMP}}.&
  \end{flalign*}
}

The next result states the linear convergence of some algorithm variants and uses the notions of pyramidal width ($\PWidth$) and minimal directional width ($\mDW$) that have been developed in \citep{lacoste2016convergence} and \citep{locatello17a} respectively, which we state in \ref{apx:basic_defs} for completeness. We note that the pyramidal width of a set $\AA$ is lower bounded by the minimal width over
all subsets of atoms, and thus is strictly greater than zero if the number of atoms is finite. The minimal directional width is a much simpler quantity and always strictly greater than zero by the symmetry of our domain.

\begin{theorem}[Linear convergence rate for strongly convex objectives] \label{thm:linear_conv_afw}
  \TheoremLinear
\end{theorem}
The previous theorem gives a geometric decrease on good steps. Combining this theorem with the bound for the number of bad steps in \eqref{eq:bad_steps}, and noting that the sufficient decrease guarantees that the objective is monotonically decreasing, we obtain a global linear convergence for \hyperref[alg:aafw]{AdaAFW}, \hyperref[alg:apfw]{AdaPFW} and \hyperref[alg:amp]{AdaMP}.

\subsection{Discussion}\label{scs:discussion}


\emph{Non-convex objectives.} \citet{lacoste2016convergence} studied the convergence of FW assuming the linear subproblems are solved exactly ($\delta=1$) and obtained a rate of the form \eqref{eq:rate_nonconvex} with $C_0 = \max\{2 h_0, L\diam(\conv(\mathcal{A}))^2\}$ instead. Both rates are similar, although our analysis is more general as it allows to consider the case in which linear subproblems are solved approximately ($\delta < 1$) and also gives rates for the Away-steps and Pairwise variants, for which no rates were previously known.

Theorem~\ref{thm:noncvx} also gives the first known convergence rates for a variant of MP on general non-convex functions. Contrary to the case of FW, this bound depends on the mean instead of the maximum of the Lipschitz estimate.

\emph{Convex objectives}.
Compared with \citep{jaggi2013revisiting}, the primal-dual rates of Theorem~\ref{thm:convex} are stronger as they hold for the last iterate and not only for the minimum over previous iterates. To the best of our knowledge, primal-dual convergence rates on the last iterate have only been derived in \citep{Nesterov2017} and were not extended to approximate linear subproblems nor the Away-steps and Pairwise variants.

Compared to \citet{Nesterov2017} on the special case of exact subproblems ($\delta=1$), the rates of Theorem~\ref{thm:convex} are similar but with $\overline{L}_t$ replaced by $L$. Hence, in the regime $\overline{L}_t \ll L$ (as is often verified in practice), our bounds have a much smaller leading constant.

For MP, \citet{locatello2018revisiting} obtain a similar convergence rate of the form $\mathcal{O}(1/(\delta^2 t))$, but with constants that depend on global properties of $\nabla f$, instead of the adaptive, averaged Lipschitz estimate in our case.

\emph{Strongly convex objectives.}
For the FW variants, the rates are identical to the ones in \citep[Theorem 1]{lacoste2015global}, with the significant difference of replacing $L$ with the adaptive~$L_t$ in the linear rate factor, giving a larger per-iteration decrease whenever $L_t < L$.
Our rates are the first also covering approximate subproblems for Away-Steps and Pairwise FW algorithms.  It's also worth noticing that both Away-steps FW and Pairwise FW have only been previously analyzed in the presence of exact line-search~\citep{lacoste2015global}.
Additionally, unlike \citep{lacoste2015global}, we do not require a smoothness assumption on $f$ outside of the domain.
Finally, for the case of MP, we again obtain the same convergence rates as in \citep[Theorem~7]{locatello17a}, but with $L$ replaced by $L_t$.

\section{Benchmarks}
\label{section:experiments}

\begin{figure*}[ht]
  \centering \begin{tabular}{lrrr|rr}
    \toprule
    {\bfseries\sffamily Dataset}                                                 & \tablefont{\#samples} & \tablefont{\#features} & {\tablefont{density}} & \multicolumn{1}{c}{$\overline{L}_{t}/L$} & \multicolumn{1}{c}{$(t - N_t)/ {t}$} \\
    \midrule
    {\bfseries\sffamily Madelon}~{\footnotesize\citep{guyon2008feature}}         & \hfill 4400           & \hfill 500             & \hfill $1.$           & \hfill $3.3 \times 10^{-3}$              & $5.0 \times 10^{-5}$                 \\
    {\bfseries\sffamily RCV1}~{\footnotesize\citep{lewis2004rcv1}}               & 697641                & \hfill 47236           & \hfill $10^{-3}$      & \hfill $1.3 \times 10^{-2}$              & $7.5 \times 10^{-5}$                 \\
    {\bfseries\sffamily MovieLens 1M}~{\footnotesize\citep{harper2016movielens}} & \hfill 6041           & \hfill 3707            & \hfill 0.04           & \hfill $1.1 \times 10^{-2}$              & --                                   \\
    \bottomrule
  \end{tabular}

  \vspace*{2em}\includegraphics[width=\linewidth]{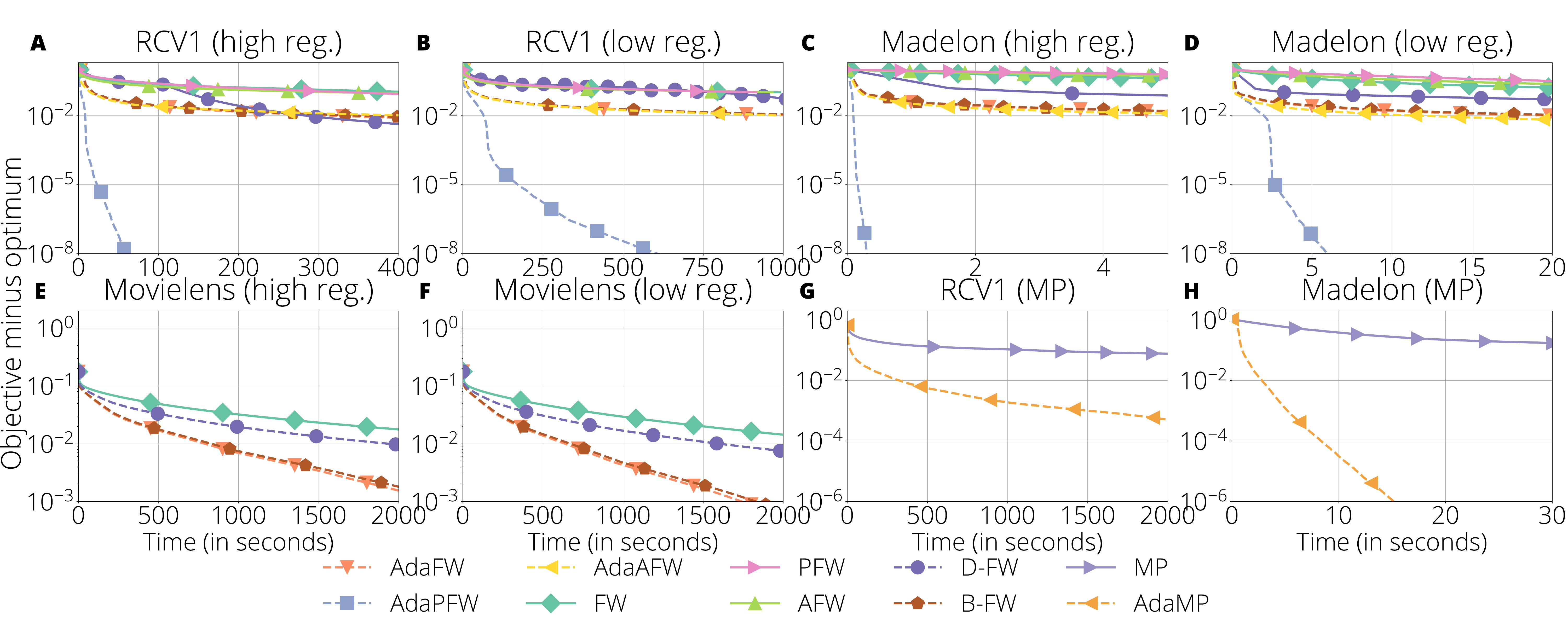}
  \caption{{\bfseries Top table}: description of the datasets. {\bfseries Bottom figure}: Benchmark of different FW and MP variants. The variants with backtracking line-search proposed in this paper are in dashed lines. Problem in A, B, C, D = logistic regression with $\ell_1$-constrained coefficients, in E, F =  Huber regression with on the nuclear norm constrained coefficients and in G, H = unconstrained logistic regression (MP variants). In all the considered datasets and regularization regimes, backtracking variants have a much faster convergence than others.}\label{fig:main_figure}
\end{figure*}

We compared the proposed methods across three problems and three datasets. The three datasets are summarized in the table of Figure~\ref{fig:main_figure}, where density denotes the fraction of nonzero coefficients in data matrix and where the last two columns are quantities that arise during the optimization of \hyperref[alg:apfw]{AdaPFW} and shed light into their empirical value. In both cases $t$ is the number of iterates until $10^{-10}$ suboptimality is achieved.

\subsection{\texorpdfstring{$\ell_1$}{l1}-constrained logistic regression}
The first problem that we consider is a logistic regression with an $\ell_1$ norm constraint on the coefficients of the form:

\begin{equation}\label{eq:l1_logistic_regression}
  \argmin_{\|\xx\|_1 \leq \beta} \ \frac{1}{n} \sum_{i=1}^n \varphi(\aa_i^\top \xx, \boldsymbol{b}_i) + \frac{\lambda}{2} \|\xx\|^2_2~,
\end{equation}
where $\varphi$ is the logistic loss. $\beta$ is chosen to give approximately $1\%$, $20\%$ of nonzero coefficients respectively. The linear subproblems in this case can be computed exactly ($\delta=1$) and consist of finding the largest entry of the gradient. The $\ell_2$ regularization parameter $\lambda$ is always set to $\lambda = \frac{1}{n}$.

We applied this problem on two different datasets: Madelon and RCV1. We show the results in Figure~\ref{fig:main_figure}, subplots A, B, C, D. In this figure we also show the performance of FW, Away-steps FW (AFW) and Pairwise FW (PFW), all of them using the step-size $ \gamma_t\!=\!\min\left\{{{g}_t}{L^{-1}\|\dd_t\|^{-2}}\!\!, \gamma_t^{\max{}}\right\}$, as well as the backtracking variants of \citet{dunn1980convergence} and \citep{beck2015cyclic}, which we denote D-FW and B-FW respectively.

\subsection{Nuclear-norm constrained Huber regression}

The second problem that we consider is collaborative filtering. We used the MovieLens 1M dataset, which contains 1 million movie ratings, and consider the problem of minimizing a Huber loss, as in  \citep{mehta2007robust}, between the true known ratings and a matrix $\boldsymbol{X}$. We also constrain the matrix by its nuclear norm $\|\boldsymbol{X}\|_* \leq \beta$, where $\beta$ is chosen to give approximately $1\%$ and $20\%$ of non-zero singular values respectively. The problem is of the form:

\begin{equation}
  \argmin_{\|\XX\|_* \leq \beta} \frac{1}{n} \sum_{(i, j) \in \mathcal{I}}^n L_\xi(\boldsymbol{A}_{i, j} -  \XX_{i, j})~,
\end{equation}
where $H_1$ is the Huber loss, defined as
\begin{equation*}
  L_\xi (a) = \begin{cases}
    \frac{1}{2}{a^2}            & \text{for } |a| \le \xi, \\
    \xi (|a| - \frac{1}{2}\xi), & \text{otherwise\,.}
  \end{cases}
\end{equation*}
The Huber loss is a quadratic for $|a|\leq \xi$ and grows linearly for $|a|>\xi$. The parameter $\xi$ controls this tradeoff and was set to 1 during the experiments.

In this case, the AFW and PFW variants were not considered as they are not directly applicable to this problem as the size of the active set is potentially unbounded. The results of this comparison can be see in subplots~E and~F of Figure~\ref{fig:main_figure}. We emphasize that the goal of this experiment is to compare different FW variants and not find the best method for matrix completion. For alternative approaches not based on FW see for instance \citep{marevcek2017matrix}.

We comment on some observed trends from these results:
\begin{itemize}[leftmargin=*]
  \item {\bfseries Importance of backtracking.} Across the different datasets, problems and regularization regimes we found that backtracking methods always perform better than their non-backtracking variant.
  \item {\bfseries Pairwise FW.} AdaPFW shows a surprisingly good performance when it is applicable, specially in the high regularization regime. A possible interpretation for this is that it is the only variant of FW in which the coefficients associated with previous atoms are not shrunk when adding a new atom, hence large step-sizes are potentially even more beneficial as coefficients that are already close to optimal do not get necessarily modified in subsequent updates.
  \item {\bfseries $\overline{L}_t$ vs $L$.} We compared the average Lipschitz estimate $\overline{L}_t$ and the $L$, the the gradient's Lipschitz constant. We found that across all datasets the former was more than an order of magnitude smaller, highlighting the need to use a local estimate of the Lipschitz constant to use a large step-size. 
  \item {\bfseries Bad steps.} Despite the very pessimistic bounds obtained for the number of bad steps in the previous section, we observe that in practice these are extremely rare events, happening less than once every 10,000 iterations.
\end{itemize}

\section{Conclusion and Future Work}
\label{section:conclusion}

In this work we have proposed and analyzed a novel adaptive step-size scheme
that can be used in projection-free methods such as FW and MP.
The method has minimal computational overhead and does not rely on any step-size  hyperparameter  (except  for  an  initial  estimate). Numerical experiments show large computational
gains on a variety of problems.

A possible extension of this work is to develop backtracking step-size strategies for randomized variants of FW such as \citep{lacoste2013block, kerdreux2018subsampling, mokhtari2018stochastic}, in which there is stochasticity in the linear subproblems.

Another area of future research is to improve the convergence rate of the Pairwise FW method. Due to the very pessimistic bound on its number of bad steps, there is still a large gap between its excellent empirical performance and its known convergence rate.
Furthermore, convergence of Pairwise and Away-steps for an infinite $\mathcal{A}$, such as the trace norm ball, is still an open problem.

\section*{Acknowledgements}

The authors would like to thank Vlad Niculae and Courtney Paquette for valuable feedback on the manuscript.

\clearpage
\bibliographystyle{plainnat}
\bibliography{biblio}

\clearpage
\appendix
\onecolumn
\titleformat{\section}{\Large\bfseries}{\thesection}{1em}{}
\titleformat{\subsection}{\large\bfseries}{\thesubsection}{1em}{}
\gdef\thesection{Appendix \Alph{section}}

{\centering{\LARGE\bfseries Linearly Convergent Frank-Wolfe with Backtracking Line-Search}

  \vspace{1em}
  \centering{{\LARGE\bfseries Supplementary material}}

}
\vspace{2em}

\paragraph{Outline.} The supplementary material of this paper is organized as follows.
\begin{itemize}[leftmargin=*]
  \item \ref{apx:pseudocode} provides pseudo-code for all FW Variants we consider: AdaFW, AdaAFW, AdaPFW, AdaMP.
  \item \ref{apx:basic_defs} contains definitions and properties relative to the objective function and/or the domain, such as the definition of geometric strong convexity and pyramidal width.
  \item \ref{apx:sufficient_decrease} we present key inequalities on the abstract algorithm which are used by the different convergence proofs.
  \item \ref{apx:proof_nonconvex} provides a proof of convergence for non-convex objectives (Theorem~\ref{thm:noncvx}).
  \item \ref{apx:proof_convex} provides a proof of convergence for convex objectives (Theorem~\ref{thm:convex}).
  \item \ref{apx:linear_proofs} provides a proof of linear convergence for all variants except FW (Theorem~\ref{thm:linear_conv_afw}).
\end{itemize}

\vspace{3em}

\section{Pseudocode}
\label{apx:pseudocode}
In this Appendix, we give detailed pseudo-code for the backtracking variants of FW (AdaFW), Away-Steps FW (AdaAFW), Pairwise FW (AdaPFW) and Matching Pursuit (AdaMP).

\subsection{Backtracking FW}

\begin{algorithm}[H]\label{alg:afw}
  \SetAlgorithmName{Algorithm}{}
  {\bfseries Input:} $\xx_0 \in \mathcal{A}$, initial Lipschitz estimate $L_{-1} > 0$, tolerance $\varepsilon \geq 0$, subproblem quality $\delta \in (0, 1]$, adaptivity params $\tau > 1, \eta \geq 1$

  \For{$t=0, 1 \ldots $}{
    Choose any $\ss_{t} \in \mathcal{A}$ that satisfies $\langle \nabla f(\xx_t),\ss_t - \xx_t\rangle \leq \delta\min_{\ss \in \mathcal{A}}\langle \nabla f(\xx_t),  \ss - \xx_t \rangle$\label{line:afw_lmo}

    Set $\dd_t=\ss_t - \xx_t$  and $\gamma_{\max} = 1$

    Set ${g}_t = \langle - \nabla f(\xx_{t}), \dd_t \rangle$

    \lIf{ ${g}_t \leq \delta\varepsilon\,$}{ \Return $\xx_t$}

    $\gamma_t, L_{t} = \text{step\_size}(f, \dd_t, \xx_t, g_t, L_{t-1}, \gamma_t^\text{max})$

    $\xx_{t+1} = \xx_t + \gamma_t \dd_t$

  }
  \caption{Backtracking FW (AdaFW)}
\end{algorithm}
\subsection{Backtracking Away-steps FW}

\begin{algorithm}[H]\label{alg:aafw}
  \SetAlgorithmName{Algorithm}{}
  {\bfseries Input:} $\xx_0 \in \mathcal{A}$, initial Lipschitz estimate $L_{-1} > 0$, tolerance $\varepsilon \geq 0$, subproblem quality $\delta \in (0, 1]$, adaptivity params $\tau > 1, \eta \geq 1$

  Let $\mathcal{S}_0 = \{\xx_0\}$ and $\balpha_{0, \vv} = 1$ for $\vv = \xx_0$ and $\balpha_{0, \vv} =0$ otherwise.

  \For{$t=0, 1 \ldots $}{
  Choose any $\ss_{t} \in \mathcal{A}$ that satisfies
  $\langle \nabla f(\xx_t),\ss_t - \xx_t\rangle \leq \delta\min_{\ss \in \mathcal{A}}\langle \nabla f(\xx_t),  \ss - \xx_t \rangle$ \label{line:aafw_lmo}

  Choose any $\vv_{t} \in \mathcal{S}_t$ that satisfies $\langle \nabla f(\xx_t),  \xx_t - \vv_t \rangle \leq \delta\min_{\vv \in {S}_t}\langle \nabla f(\xx_t), \xx_t - \vv\rangle$ \label{line:away_lmo}

  \eIf{$\langle \nabla f(\xx_t) , \ss_t - \xx_t \rangle \leq \langle \nabla f(\xx_t) ,\xx_t - \vv_t\rangle$ \label{line:decision_step}}{

  $\dd_t=\ss_t - \xx_t$ and $\gamma_t^{\text{max}}=1$\label{l:gammamax}
  }{
  $\dd_t=\xx_t - \vv_t$,  and $\gamma_t^{\text{max}}\!=\!{\alpha_{\vv_t, t}}/{(1\!-\!\alpha_{\vv_t, t})}$\label{l:gammamax2}
  }

  Set ${g}_t = \langle - \nabla f(\xx_{t}), \dd_t \rangle$

  \lIf{ ${g}_t \leq \delta\varepsilon\,$}{ \Return $\xx_t$}

  $\gamma_t, L_{t} = \text{step\_size}(f, \dd_t, \xx_t, g_t, L_{t-1}, \gamma_t^\text{max})$

  $\xx_{t+1} = \xx_t + \gamma_t \dd_t$

  Update active set $\mathcal{S}_{t+1}$ and $\balpha_{t+1}$ (see text)
  }
  \caption{Backtracking Away-Steps FW (AdaAFW)}
\end{algorithm}
The active set is updated as follows.

\begin{itemize}
  \item In the case of a FW step, we update the support set $\mathcal{S}_{t+1}=\{\ss_t\}$ if $\gamma_t=1$ and otherwise $\mathcal{S}_{t+1}=\mathcal{S}_{t}\cup\{\ss_t\}$, with coefficients $\alpha_{\vv, t+1}=(1-\gamma_t)\alpha_{\vv, t}$ for $\vv\in\mathcal{S}_t\setminus\{\ss_t\}$ and $\alpha_{\ss_t, t+1}=(1-\gamma_t)\alpha_{\ss_t, t}+\gamma_t$.
  \item In the case of an Away step:
        If $\gamma_t=\gamma_{\text{max}}$, then
        $\mathcal{S}_{t+1}=\mathcal{S}_{t}\setminus\{\vv_t\}$, and if $\gamma_t < \gamma_{\text{max}}$,
        then  $\mathcal{S}_{t+1}=\mathcal{S}_{t}$. Finally, we update the weights as $\alpha_{\vv, t+1}={(1+\gamma_t)\alpha_{\vv, t}}$ for $\vv\in\mathcal{S}_{t}\setminus \{\vv_t\}$ and $\alpha_{\vv_t, t+1}={(1+\gamma_t)}\alpha_{\vv_t, t}-\gamma_t$ for the other atoms.
\end{itemize}

\subsection{Backtracking Pairwise FW}

\begin{algorithm}[H]\label{alg:apfw}
  {\bfseries Input:} $\xx_0 \in \mathcal{A}$, initial Lipschitz estimate $L_{-1} > 0$, tolerance $\varepsilon \geq 0$, subproblem quality $\delta \in (0, 1]$, adaptivity params $\tau > 1, \eta \geq 1$

  Let $\mathcal{S}_0 = \{\xx_0\}$ and $\balpha_{0, \vv} = 1$ for $\vv = \xx_0$ and $\balpha_{0, \vv} =0$ otherwise.

  \For{$t=0, 1 \ldots $}{
  Choose any $\ss_{t} \in \mathcal{A}$ that satisfies
  $\langle \nabla f(\xx_t),\ss_t - \xx_t\rangle \leq \delta\min_{\ss \in \mathcal{A}}\langle \nabla f(\xx_t),  \ss - \xx_t \rangle$\label{line:apfw_lmo}

  Choose any $\vv_{t} \in \mathcal{S}_t$ that satisfies $\langle \nabla f(\xx_t),  \xx_t - \vv_t \rangle \leq \delta\min_{\vv \in {S}_t}\langle \nabla f(\xx_t), \xx_t - \vv\rangle$ \label{line:plmo}

  $\dd_t=\ss_t - \vv_t$ and $\gamma_t^{\text{max}}\!=\!{\alpha_{\vv_t, t}}$

  Set ${g}_t = \langle - \nabla f(\xx_{t}), \dd_t \rangle$

  \lIf{ ${g}_t \leq \delta\varepsilon\,$}{ \Return $\xx_t$}

  $\gamma_t, L_{t} = \text{step\_size}(f, \dd_t, \xx_t, g_t, L_{t-1}, \gamma_t^\text{max})$

  $\xx_{t+1} = \xx_t + \gamma_t \dd_t$

  Update active set $\mathcal{S}_{t+1}$ and $\balpha_{t+1}$ (see text)
  }
  \caption{Backtracking Pairwise FW (AdaPFW)}
\end{algorithm}

\hyperref[alg:apfw]{AdaPFW} only moves weight from $\vv_t$ to $\ss_t$. The active set update becomes $\alpha_{\ss_t, t+1}=\alpha_{\ss_t, t}+\gamma_t$, $\alpha_{\vv_t, t+1}=\alpha_{\vv_t, t}-\gamma_t$, with $\mathcal{S}_{t+1}=(\mathcal{S}_{t}\setminus \{\vv_t\})\cup\{\ss_t\}$ if $\alpha_{\vv_t, t+1} = 0$ and $\mathcal{S}_{t+1}=\mathcal{S}_{t}\cup\{\ss_t\}$ otherwise.

\subsection{Backtracking Matching Pursuit}\label{apx:pseudocode_mp}

Matching Pursuit~\citep{mallat1993matching, locatello17a} is an algorithm to solve optimization problems of the form
\begin{empheq}[box=\mybluebox]{equation}\label{eq:obj_fun_mp}
  \minimize_{\xx \in \lin({\mathcal{A}})}\,\vphantom{\sum_i^n} f(\xx) ~,
\end{empheq}
where $\lin(\AA) \defas \left \{ {\sum_{\vv \in \mathcal{A}} \lambda_\vv \vv \,\big|\, \lambda_\vv  \in \RR} \right \}$ is the linear span of the set of  \emph{atoms} $\mathcal{A}$.
As for the backtracking FW algorithm, we assume that $f$ is $L$-smooth and $\mathcal{A}$ a potentially infinite but bounded set of elements in $\RR^p$.

The MP algorithm relies on solving at each iteration a linear subproblem over the set $\mathcal{B} \defas\mathcal{A}\cup -\AA$, with $-\AA = \{-\aa \,|\, \aa \in \mathcal{A}\}$. The linear subproblem that needs to be solved at each iteration is the following, where as for previous variants, we allow for an optional quality parameter $\delta \in (0, 1]$:
\begin{equation}\label{eq:mp_lmo}
\langle \nabla f(\xx_t), \ss_t\rangle \leq \delta\min_{\ss \in \mathcal{B}}\langle \nabla f(\xx_t), \ss\rangle~.
\end{equation}
In Algorithm~\ref{alg:amp} we detail a novel adaptive variant of the MP algorithm, which we name \hyperref[alg:amp]{AdaMP}.

\begin{algorithm}[H]\label{alg:amp}
  {\bfseries Input:} $\xx_0 \in \mathcal{A}$, initial Lipschitz estimate $L_{-1} > 0$, tolerance $\varepsilon \geq 0$, subproblem quality $\delta \in (0, 1]$

  \For{$t=0, 1 \ldots $}{
  Choose any $\ss_{t} \in \mathcal{A}$ that satisfies \eqref{eq:mp_lmo}\label{line:amp_lmo}

  $\dd_t=\ss_t$

  Set ${g}_t = \langle - \nabla f(\xx_{t}), \dd_t \rangle$

  \lIf{ ${g}_t \leq \delta\varepsilon\,$}{ \Return $\xx_t$}

  $\gamma_t, L_{t} = \text{step\_size}(f, \dd_t, \xx_t, g_t, L_{t-1}, \infty)$

  $\xx_{t+1} = \xx_t + \gamma_t \dd_t$

  }
  \caption{Backtracking Matching Pursuit (AdaMP)}
\end{algorithm}





\clearpage

\section{Basic definitions and properties}\label{apx:basic_defs}

In this section we give basic definitions and properties relative to the objective function and/or the domain, such as the definition of geometric strong convexity and pyramidal width.
These definitions are not specific to our algorithms and have appeared in different sources such as \citet{lacoste2015global, locatello17a}. We merely gather them here for completeness.



\begin{definition}[Geometric strong convexity] We define the \textbf{geometric strong convexity constant $\mu_f^A$} as
  \begin{align}
     & \mu_f^{\mathcal{A}} \defas \inf_{\substack{\xx, \xx^\star \in \conv(\mathcal{A})                                                                                                        \\ \langle \nabla f(\xx), \xx^\star - \xx \rangle < 0}} \dfrac{2}{\gamma(\xx,\xx^\star)^2} \Big( f(\xx^\star) - f(\xx) - \langle \nabla f(\xx), \xx^\star - \xx \rangle \Big) \label{eq:geo-strong-cvx} \\
     & \text{ where }\gamma (\xx,\xx^\star) \defas \dfrac{ \langle -\nabla f(\xx), \xx^\star - \xx \rangle}{\langle -\nabla f(\xx), \ss_f(\xx) - \vv_f(\xx) \rangle }, \;\;
  \end{align}
  where
  \begin{align}
    \ss_f(\xx)           & \defas \argmin_{\vv \in \mathcal{A}} \langle \nabla f(\xx), \vv \rangle \\
    \vv_f(\xx)           & \defas \argmin_{\substack{\vv = \vv_\mathcal{S}(\xx)                    \\ \mathcal{S} \in \mathcal{S}_{\xx}}} \langle \nabla f(\xx), \vv \rangle \\
    \vv_\mathcal{S}(\xx) & \defas \argmax_{\vv \in \mathcal{S}} \langle \nabla f(\xx), \vv \rangle
  \end{align}
  where $\mathcal{S} \subseteq \AA$ and $\mathcal{S}_{\xx} \defas \{\mathcal{S} | \mathcal{S} \subseteq \mathcal{A}$ such that $\xx$ is a proper convex combination of all the elements in $\mathcal{S}\}$ (recall $\xx$ is a proper convex combination of elements in $\SS$ when $\xx = \sum_i \alpha_i \ss_i$ where $\ss_i \in \SS$ and $\alpha_i \in (0,1)$).
\end{definition}

\begin{definition}[Pyramidal width] The \textbf{pyramidal width} of a set $\AA$ is the smallest pyramidal width of all its faces, i.e.
  \begin{align}
    \PWidth(\AA) \defas \min_{\substack{ \xx \in \mathcal{K} \\ \mathcal{K} \in \face(\conv(\AA)) \\ \rr \in \cone(\mathcal{K} - \xx) \backslash \{ 0 \}}} \PdirW(\mathcal{K} \cap \mathcal{A}, \rr,\xx)
  \end{align}
  where $\PdirW$ is the pyramidal directional width, defined as
  \begin{align}
    \PdirW(W)(\AA,\rr,\xx) \defas \min_{\mathcal{S} \in \mathcal{S}_x} \max_{\ss \in \AA, \vv \in \mathcal{S}} \Big \langle \dfrac{\rr}{\| \rr\|_2}, \ss - \vv \Big \rangle
  \end{align}
\end{definition}

We now relate these two geometric quantities together.
\begin{lemma}[Lower bounding $\mu_f^A$] \label{lemma:lb_mufa} Let $f$ $\mu$--strongly convex on $\conv(\mathcal{A}) = \conv (\AA)$. Then
  \begin{align}
    \mu_f^A \geq \mu \cdot (\PWidth(\AA))^2
  \end{align}
\end{lemma}

\begin{proof}
  We refer to \citep[Theorem 6]{lacoste2015global}.
\end{proof}

\begin{proposition} \label{prop:ub:pwdith} $\PWidth(\AA ) \leq \diam(\conv(\AA))$ where $\diam(\mathcal{X}) \defas \sup_{x,y \in \mathcal{X}} \|x - y\|_2$.
\end{proposition}

\begin{proof}
  First note that given $\rr \in \mathcal{R}$, $\ss \in \mathcal{S}$, $\vv \in \mathcal{V}$ with $\mathcal{R}$, $ \mathcal{S}$, $\mathcal{V} \subseteq \mathbb{R}^n$, we have
  \begin{align}
    \langle \rr/\|\rr\|_2, \ss - \vv \rangle                                                                                        & \leq \|\ss - \vv\|_2                                               & \forall \rr \in \mathcal{R},\ss \in \mathcal{S},\vv \in \mathcal{V} \\
    \Rightarrow  \max_{\ss \in \mathcal{S}, \vv \in \mathcal{V}} \langle \rr/\|\rr\|_2, \ss - \vv \rangle                           & \leq \max_{\ss \in \mathcal{S}, \vv \in \mathcal{V}} \|\ss-\vv\|_2 & \forall \rr \in \mathcal{R}                                         \\
    \Rightarrow \min_{\rr \in \mathcal{R}} \max_{\ss \in \mathcal{S}, \vv \in \mathcal{V}} \langle \rr/\|\rr\|_2, \ss - \vv \rangle & \leq \max_{\ss \in \mathcal{S}, \vv \in \mathcal{V}} \|\ss-\vv\|_2
  \end{align}

  Applying this result to the definition of pyramidal width we have
  \begin{align}
    \PWidth(\AA) & = \min_{\substack{ \xx \in \mathcal{K}                                                                                    \\ \mathcal{K} \in \face(\conv(\AA)) \\ \rr \in \cone(\mathcal{K} - \xx) \backslash \{ 0 \}}}\PdirW(\mathcal{K} \cap \mathcal{A}, \rr,\xx) \\
                 & = \min_{\substack{ \xx \in \mathcal{K}                                                                                    \\ \mathcal{K} \in \face(\conv(\AA)) \\ \rr \in \cone(\mathcal{K} - \xx) \backslash \{ 0 \}}} \min_{\mathcal{S} \in \mathcal{S}_x} \max_{\ss \in \AA, \vv \in \mathcal{S}} \Big \langle \dfrac{\rr}{\| \rr\|}, \ss - \vv \Big \rangle \\
                 & = \min_{r \in \mathcal{R}} \max_{s \in \AA, v \in \mathcal{V}} \Big \langle \dfrac{\rr}{\| \rr\|}, \ss - \vv \Big \rangle \\
  \end{align}
  where $\mathcal{R} = \{ \cone (\mathcal{K} - \xx) \backslash \{0\} : \text{for some } \xx \in \mathcal{K}, \;  \mathcal{K} \in \face(\conv(\AA))\}$ and $\mathcal{V}$ is some subset of $\AA$. Applying the derived result we have that
  \begin{align*}
    \PWidth(\AA) & \leq \max_{\ss \in \AA, \vv \in \mathcal{V}} \| \ss - \vv\|_2 \\
                 & \leq \max_{\ss, \vv \in \conv(\AA)} \| \ss - \vv\|_2          \\
                 & = \diam(\conv(\AA))
  \end{align*}
\end{proof}

\begin{definition} The \textbf{minimal directional width} $\mDW(\AA)$ of a set of atoms $\AA$ is defined as
  \begin{align}
    \mDW (\AA) = \min_{\dd \in \lin(\AA)} \max_{\zz \in \AA} \dfrac{\langle \zz,\dd \rangle}{\|\dd\|}~.
  \end{align}
\end{definition}

Note that in contrast to the pyramidal width, the minimal directional width here is a much simpler and robust property of the atom set $\AA$, not depending on its combinatorial face structure of the polytope. As can be seen directly from the definition above, the $\mDW(\AA)$ is robust when adding a duplicate atom or small perturbation of it to $\AA$.

\clearpage

\section{Preliminaries: Key Inequalities}\label{apx:sufficient_decrease}

In this appendix we prove that the sufficient decrease condition verifies a recursive inequality. This key result is used by all convergence proofs.

\begin{lemma}\label{lemma:sufficient_decrease_inequality}
  The following inequality is verified for all proposed algorithms (with $\gamma_t^{\max} = +\infty$ for \hyperref[alg:amp]{AdaMP}):
  \begin{equation}
    f(\xx_{t+1}) \leq f(\xx_t)  -\xi {g_t}+ \frac{\xi^2 L_{t}}{2}\|\dd_t\|^2~ \text{ for all }\xi \in [0, \gamma_t^{\max}].
  \end{equation}
\end{lemma}

\begin{proof}
  We start the proof by proving an optimality condition of the step-size. Consider the following quadratic optimization problem:
  \begin{equation}\label{eq:step_size_optimality}
    \minimize_{\xi \in [0, \gamma_t^{\max}]} -\xi {g_t} + \frac{{L}_t \xi^2}{2} \|\dd_t\|^2 ~.
  \end{equation}
  Deriving with respect to $\xi$ and noting that on all the considered algorithms we have $\langle \nabla f(\xx_t), \dd_t \rangle \leq 0$, one can easily verify that the global minimizer is achieved at the value
  \begin{equation}
    \displaystyle\min\left\{\frac{{g_t}}{L_t\|\dd_t\|^{2} }, \gamma_t^{\text{max}}\right\}~,
  \end{equation}
  where ${g_t} = \langle - \nabla f(\xx), \dd_t\rangle$.
  This coincides with the value of $\gamma_{t+1}$ computed by the backtracking procedure on the different algorithms and so we have:
  \begin{equation}
    -\gamma_t{g_t} + \frac{{L}_t {\gamma_t}^2}{2} \|\dd_t\|^2  \leq  -\xi {g_t} + \frac{{L}_t \xi^2}{2} \|\dd_t\|^2   ~\text{ for all $\xi \in [0, \gamma^{\max}]$}~.
  \end{equation}

  We can now write the following sequence of inequalities, that combines the sufficient decrease condition with this last inequality:
  \begin{align}
    f(\xx_{t+1}) & \,~{\leq}~\, f(\xx_t) -\gamma_t {g_t} + \frac{{L}_t \gamma_t^2}{2}\|\dd_t\|^2 \quad                                                                                                   \\
                 & \stackrel{\eqref{eq:step_size_optimality}}{\leq} f(\xx_t) - \xi {g_t} + \frac{{L}_t \xi^2}{2}\|\dd_t\|^2\label{eq:recusive_rhs_final}~ \text{ for any $\xi \in [0, \gamma^{\max}]$}~.
  \end{align}
\end{proof}

\begin{proposition}\label{prop:lipschitz_estimate_bound}
  The Lipschitz estimate ${L}_t$ is bounded as ${L}_t \leq \max\{\tau L, L_{-1}\}$.
\end{proposition}
\begin{proof}

  If the sufficient decrease condition is verified then we have $L_t = \eta L_{t-1}$ and so $L_t \leq L_{t-1}$. If it's not, we at least have that the Lipschitz estimate cannot larger than $\tau L$ by definition of Lipschitz constant. Combining both bounds we obtain
  \begin{equation}
    L_t \leq \max\{\tau L, L_{t-1}\}~.
  \end{equation}
  Applying the same bound recursively on $L_{t-1}$ leads to the claimed bound $L_t \leq \max\{\tau L, L_{-1}\}$.
\end{proof}

\begin{lemma}\label{lemma:bound_gt_fw_gap}
  Let $g(\cdot)$ be as in Theorem~\ref{thm:noncvx}, i.e.,  $g(\cdot) = g^{\text{FW}}(\cdot)$ for FW variants (\hyperref[alg:afw]{AdaFW}, \hyperref[alg:aafw]{AdaAFW}, \hyperref[alg:apfw]{AdaPFW}) and  $g(\cdot) = g^{\text{MP}}(\cdot)$ for MP variants (\hyperref[alg:amp]{AdaMP}). Then for any of these algorithms we have
  \begin{equation}\label{eq:bound_fw_gap_gt}
    {g}_t  \geq \delta g(\xx_t)~.
  \end{equation}
\end{lemma}
\begin{proof}

  \begin{itemize}[leftmargin=*]
    \item For \hyperref[alg:afw]{AdaFW} and \hyperref[alg:amp]{AdaMP}, Eq. \eqref{eq:bound_fw_gap_gt} follows immediately from the definition of ${g_t}$ and $g(\xx_t)$.
    \item For \hyperref[alg:aafw]{AdaAFW}, by the way the descent direction is selected in Line \ref{line:decision_step}, we always have
          \begin{equation}
            {g}_t \geq \langle \nabla f(\xx_t), \xx_t - \ss_t \rangle \geq \delta g(\xx_t)~,
          \end{equation}
          where the last inequality follows from the definition of $\ss_t$
    \item For \hyperref[alg:apfw]{AdaPFW}, we have
          \begin{align}
            {g}_t & = \langle \nabla f(\xx_t), \vv_t - \ss_t \rangle =  \langle \nabla f(\xx_t), \xx_t - \ss_t \rangle  +  \langle \nabla f(\xx_t), \vv_t - \xx_t \rangle \\
                  & \geq \langle \nabla f(\xx_t), \xx_t - \ss_t \rangle \geq \delta g(\xx_t)
          \end{align}
          where the term $\langle \nabla f(\xx_t), \vv_t - \xx_t \rangle$ is positive by definition of $\vv_t$ since $\xx_t$ is necessarily in the convex envelope of $\mathcal{S}_t$. The second inequality follows from the definition of $\ss_t$.
  \end{itemize}
\end{proof}

\begin{customtheorem}{\ref{thm:bound_backtracking}}
  Let $N_t$ be the total number of evaluations of the sufficient decrease condition up to iteration $t$. Then we have
  \begin{equation}
    n_t \leq \left[ 1 -\frac{\log \eta}{\log \tau}\right](t+1) + \frac{1}{\log\tau}\max\left\{\log\frac{\tau L}{L_{-1}}, 0\right\}~.
  \end{equation}
\end{customtheorem}
\begin{proof}
  This proof follows roughly that of \citep[Lemma 3]{nesterov2013gradient}, albeit with a slightly different bound on $L_t$ due to algorithmic differences.

  Denote by $n_i \geq 1$ the number of evaluations of the sufficient decrease condition. Since the algorithm multiplies by $\tau$ every time that the sufficient condition is not verified, we have
  \begin{equation}
    L_{i} = \eta L_{i-1} \tau^{n_i-1} ~.
  \end{equation}
  Taking logarithms on both sides we obtain
  \begin{equation}
    n_i \leq 1  - \frac{\log \eta}{\log \tau} + \frac{1}{\tau} \log\frac{L_{i}}{L_{i-1}}~.
  \end{equation}
  Summing from $i=0$ to $i=t$ gives
  \begin{equation}
    n_t \leq \sum_{i=0}^t n_i = \left[ 1 -\frac{\log \eta}{\log \tau}\right](t+1) + \frac{1}{\log\tau}\log\left(\frac{L_t}{L_{-1}}\right)
  \end{equation}
  Finally, from Proposition \ref{prop:lipschitz_estimate_bound} we have the bound $L_t \leq \max\{\tau L, L_{-1}\}$, which we can use to bound the numerator's last term. This gives the claimed bound
  \begin{equation}
    n_t \leq \sum_{i=0}^t n_i = \left[ 1 -\frac{\log \eta}{\log \tau}\right](t+1) + \frac{1}{\log\tau}\max\left\{\log\frac{\tau L}{L_{-1}}, 0\right\}~.
  \end{equation}
\end{proof}

\subsection{A bound on the number of bad steps}\label{scs:bound_bad_steps}

To prove the linear rates for the backtracking AFW and backtracking PFW algorithm it is necessary to bound the number of bad steps. There are two different types of bad steps: ``drop'' steps and ``swap'' steps. These names come from how the active set $\SS_t$ changes. In a drop step, an atom is removed from the active set (i.e. $|\SS_{t+1}| < |\SS_{t}|$). In a swap step, the size of the active set remains unchanged (i.e. $|\mathcal{S}_{t+1}| = |\mathcal{S}_{t}|$) but one atom is swapped with another one not in the active set. Note that drop steps can occur in the (backtracking) Away-steps and Pairwise, but swap steps can only occur in the Pairwise variant.

For the proofs of linear convergence in \ref{apx:linear_proofs}, we show that these two types of bad steps are only problematic when $\gamma_t = \gamma_t^{\max} < 1$. In these scenarios, we cannot provide a meaningful decrease bound. However, we show that the number of bad steps we take is bounded. The following two lemmas adopted from \citep[Appendix C]{lacoste2015global} bound the number of drop steps and swap steps the backtracking algorithms can take.

\begin{lemma}\label{lemma:drop_steps}
  After $T$ steps of \hyperref[alg:aafw]{AdaAFW} or \hyperref[alg:apfw]{AdaPFW}, there can only be ${T}/{2}$ drop steps. Also, if there is a drop step at step $t+1$, then $f(\xx_{t+1})- f(\xx_{t}) < 0$.
\end{lemma}

\begin{proof}
  Let $A_t$ denote the number of steps that added a vertex in the expansion, and let $D_t$ be the number of drop steps. Then $1 \leq |\SS_{t}| = |\SS_{0}| + A_t - D_t$ and we clearly have $A_t - D_t \leq t$. Combining these two inequalities we have that $D_t \leq \frac{1}{2} (|\SS_{0}| - 1 + t) = \frac{t}{2}$.

  To show $f(\xx_{t+1})- f(\xx_{t}) < 0$, because of Lemma~\ref{lemma:sufficient_decrease_inequality}, it suffices to show that
  \begin{align}
    -\gamma_t g_t + \dfrac{1}{2} \gamma_t^2 L_t \|\dd_t\|^2 < 0~,
  \end{align}
  with $\gamma_t = \gamma_t^{\max}$ (recall drop steps only occur when $\gamma_t = \gamma_t^{\max}$). Note this is a convex quadratic in $\gamma_t$ which is precisely less than or equal to $0$ when $\gamma_t \in [0,2g_t /L_t\|\dd_t\|^2]$. Thus in order to show $f(\xx_{t+1})- f(\xx_{t}) < 0$ it suffices to show $\gamma_t^{\max} \in (0,2g_t /L_t\|\dd_t\|^2)$. This follows immediately since $0 < \gamma_t^{\max} \leq g_t/L_t \|\dd_t\|^2$.
\end{proof}

Since in the \hyperref[alg:aafw]{AdaAFW} algorithm all bad steps are drop steps, the previous lemma implies that we can effectively bound the number of bad steps by $t/2$, which is the bound claimed in \eqref{eq:bad_steps}.

\begin{lemma}\label{lemma:swap_steps}
  There are at most $3|\AA|!$ bad steps between any two good steps in \hyperref[alg:apfw]{AdaPFW}. Also, if there is a swap step at step $t+1$, then $f(\xx_{t+1})- f(\xx_{t}) < 0$.
\end{lemma}

\begin{proof}
  Note that bad steps only occur when $\gamma_t = \gamma_t^{\max} = \alpha_{\vv_t,t}$. When this happens there are two possibilities; we either move all the mass from $\vv_t$ to a new atom $\ss_t \not \in \SS_t$ (i.e. $\alpha_{\vv_t, t+1} = 0$ and $\alpha_{\ss_t,t+1} = \alpha_{\vv_t,t}$ ) and preserve the cardinality of our active set ($|\SS_{t+1}| = |\SS_t|$) or we move all the mass from $\vv_t$ to an old atom $\ss_t \in \SS_t$ (i.e. $\alpha_{\ss_t,t+1} = \alpha_{\ss_t,t} + \alpha_{\vv_t,t}$) and the cardinality of our active set decreases by 1 ($|\SS_{t+1}| < |\SS_t|$).
  In the former case, the possible values of the coordinates $\alpha_{\vv}$ do not change, but they are simply rearranged in the possible $|\AA|$ slots. Note further every time the mass from $\vv_t$ moves to a new atom $\ss_t \not \in \SS_t$ we have strict descent, i.e. $f(\xx_{t+1}) < f(\xx_t)$ unless $\xx_t$ is already optimal (see Lemma \ref{lemma:drop_steps}) and hence we cannot revisit the same point unless we have converged. Thus the maximum number of possible consecutive swap steps is bounded by the number of ways we can assign $|\SS_t|$ numbers in $|\AA|$ slots, which is $|\AA|! /(|\AA| - |\SS_t|)!$. Furthermore, when the cardinality of our active set drops, in the worst case we will do a maximum number of drop steps before reducing the cardinality of our active set again. Thus starting with $|\SS_t| = r$ the maximum number of bad steps $B$ without making any good steps is upper bounded by
  \begin{align*}
    B \leq \sum_{k=1}^r \dfrac{|\AA|!}{(|\AA| - k)!} \leq |\AA|! \sum_{k=0}^\infty \dfrac{1}{k!} = |\AA|! e \leq 3|\AA|!
  \end{align*}
\end{proof}



\clearpage

\section{Proofs of convergence for non-convex objectives}\label{apx:proof_nonconvex}

In this appendix we provide the convergence proof of Theorem~\ref{thm:noncvx}. Although this theorem provides a unified convergence proof for both variants of FW and MP, for convenience we split the proof into one for FW variants (Theorem \ref{thm:noncvx}.A) and another one for variants of MP (Theorem \ref{thm:noncvx}.B)
\vspace{1em}

\begin{mdframed}
  \begin{customtheorem}{\ref{thm:noncvx}.A}
    Let $\xx_{t}$ denote the iterate generated by either \hyperref[alg:afw]{AdaFW}, \hyperref[alg:aafw]{AdaAFW} or \hyperref[alg:apfw]{AdaPFW} after $t$ iterations. Then for any iteration $t$ with $N_{t+1}\geq 0$, we have the following suboptimality bound in terms of the FW gap:
    \begin{equation}
      \lim_{k\to \infty}g^{\text{FW}}(\xx_k) = 0\quad\text{ and }\quad \min_{k=0,\dots, t} g^{\text{FW}}(\xx_k) \leq \frac{\max\{2 h_0, {L^{\max}_t}\diam(\mathcal{A})^2\}}{\delta\sqrt{N_{t+1}}} = \mathcal{O}\left(\frac{1}{\delta\sqrt{t}}\right)
    \end{equation}
  \end{customtheorem}
\end{mdframed}
\begin{proof}
  By Lemma~\ref{lemma:sufficient_decrease_inequality} we have the following inequality for any $k$ and any $\xi \in [0,\gamma_k^{\max{}}]$,
  \begin{align}\label{eq:nonconvex_sfd}
    f(\xx_{k+1}) & \leq f(\xx_k) - \xi {g}_k + \frac{\xi^2 C_k}{2}~,
  \end{align}
  where we define $C_k \defas L_k\|\dd_k\|^2$ for convenience.
  We consider now different cases according to the relative values of $\gamma_k$ and $\gamma_k^{\max}$, yielding different upper bounds for the right hand side.

  \textbf{\underline{Case 1: $\gamma_k < \gamma_k^{\max}$}}\\
  In this case, $\gamma_k$ maximizes the right hand side of the (unconstrained) quadratic in inequality \eqref{eq:nonconvex_sfd} which then becomes:

  \begin{align}
    f(\xx_{k+1}) & \leq f(\xx_k) -\frac{{g}_k^2}{2C_k} \leq f(\xx_k) - \frac{{g}_k}{2}\min\left\{\frac{{g}_k}{C_k}, 1\right\}
  \end{align}

  \textbf{\underline{Case 2: $\gamma_k =\gamma_k^{\max} \geq 1$}} \\
  By the definition of $\gamma_t$, this case implies that $C_k \leq g_k$ and so using $\xi=1$ in \eqref{eq:nonconvex_sfd} gives
  \begin{equation}
    f(\xx_{k+1}) - f(\xx_k) \leq -{g}_k + \frac{C_k}{2} \leq -\frac{{g}_k}{2}~.
  \end{equation}


  \textbf{\underline{Case 3: $\gamma_k =\gamma_k^{\max} < 1$}}

  This corresponds to the problematic drop steps for \hyperref[alg:aafw]{AdaAFW} or possibly swap steps for \hyperref[alg:apfw]{AdaPFW}, in which we will only be able to guarantee that the iterates are non-increasing. Choosing $\xi=0$ in \eqref{eq:nonconvex_sfd} we can at least guarantee that the objective function is non-increasing:
  \begin{equation}
    f(\xx_{k+1}) - f(\xx_k) \leq 0~.
  \end{equation}

  \paragraph{Combining the previous cases.}
  We can combine the inequalities obtained for the previous cases into the following inequality, valid for all $k\leq t$,
  \begin{align}
    \label{ineq:AFW_NC_obj_decrease}
    f(\xx_{k+1}) - f(\xx_k) & \leq  - \frac{{g}_k}{2}\min\left\{\frac{{g}_k}{C_k}, 1\right\} \mathds{1}\{k \text{ is a good step} \}
  \end{align}
  Adding the previous inequality from $k=0$ up to $t$ and rearranging we obtain
  \begin{align}
    f(\xx_0) - f(\xx_{t+1}) & \geq \sum_{k=0}^t \frac{{g}_k}{2}\min\left\{\frac{{g}_k}{L_k \|\dd_k\|^2}, 1\right\}\mathds{1}\{k \text{ is a good step} \}                \\
                            & \geq \sum_{k=0}^t \frac{{g}_k}{2}\min\left\{\frac{{g}_k}{{C}_k^{\max}}, 1\right\}\mathds{1}\{k \text{ is a good step} \}\label{eq:sum_g_k}
  \end{align}
  with $C_t^{\max} \defas L_t^{\max} \diam(\conv(\mathcal{A}))^2$.
  Taking the limit for $t \to +\infty$ we obtain that the left hand side is bounded by the compactness assumption on the domain $\conv(\mathcal{A})$ and $L$-smoothness on $f$. The right hand side is an infinite sum, and so a necessary condition for it to be bounded is that $g_k \to 0$, since $g_k\geq 0$ for all $k$. We have hence proven that $\lim_{k\to \infty} g_k = 0$, which by Lemma \ref{lemma:bound_gt_fw_gap} implies $\lim_{k\to \infty} g(\xx_k) = 0$. This proves the first claim of the Theorem.

  We will now aim to derive explicit convergence rates for convergence towards a stationary point. Let $\widetilde{g_t} = \min_{0 \leq k \leq t}{{g}_k}$, then from Eq.~\eqref{eq:sum_g_k} we have
  \begin{align}
    f(\xx_0) - f(\xx_{t+1}) & \geq \sum_{k=0}^t \frac{\widetilde{g_t}}{2}\min\left\{\frac{\widetilde{g_t}}{{C}_t^{\max}}, 1\right\}\mathds{1}\{k \text{ is a good step} \} \\
                            & = N_{t+1} \frac{{\widetilde{g_t}}}{2}\min\left\{\frac{\widetilde{g_t}}{{C}_t^{\max}}, 1\right\}
    ~.\label{eq:nonconvexsum}
  \end{align}

  We now make a distinction of cases for the quantities inside the min.
  \begin{itemize}[leftmargin=*]
    \item If $\widetilde{g_t}\leq {C}_t^{\max}$, then \eqref{eq:nonconvexsum} gives $f(\xx_0) - f(\xx_{t+1}) \geq N_{t+1} {\widetilde{g_t}^2}/{(2 {C}_t^{\max})}$, which reordering  gives
          \begin{align}
            \widetilde{g_t}\leq \sqrt{\frac{2{C}_t^{\max} (f(\xx_0) - f(\xx_{t+1}))}{N_{t+1}}} \leq \sqrt{\frac{2{C}_t^{\max} h_0}{N_{t+1}}}
            \leq \frac{2 h_0 +  {C}_t^{\max}}{2\sqrt{N_{t+1}}}  \leq \frac{\max\{2 h_0, {C}_t^{\max}\}}{\sqrt{N_{t+1}}} ~.
          \end{align}
          where in the third inequality we have used the inequality $\sqrt{ab}\leq \frac{a + b}{2}$ with $a = \sqrt{2 h_0}$, $b = \sqrt{{C}_t^{\max}}$.
    \item If $\widetilde{g_t}>{C}_t^{\max}$ we can get a better $\frac{1}{N_t}$ rate, trivially bounded by $\frac{1}{\sqrt{N_t}}$.
          \begin{align}
            \widetilde{g_t} \leq \frac{2  h_0}{N_{t+1}} \leq \frac{2 h_0}{\sqrt{N_{t+1}}} \leq \frac{\max\{2 h_0, {C}_t^{\max}\}}{\sqrt{N_{t+1}}} ~.
          \end{align}

  \end{itemize}

  We have obtained the same bound in both cases, hence we always have
  \begin{align}
    \widetilde{g_t} \leq  \frac{\max\{2 h_0, {C}_t^{\max}\}}{\sqrt{N_{t+1}}}~.
  \end{align}
  Finally, from Lemma~\ref{lemma:bound_gt_fw_gap} we have $g(\xx_k) \leq \frac{1}{\delta} g_k$ for all $k$ and so
  \begin{equation}
    \min_{0 \leq k \leq t} g(\xx_k) \leq \frac{1}{\delta}\min_{0 \leq k \leq t} g_k = \frac{1}{\delta} \widetilde{g_t} \leq \frac{\max\{2 h_0, {C}_t^{\max}\}}{\delta\sqrt{N_{t+1}}} ~,
  \end{equation}
  and the claimed bound follows by definition of ${C}_t^{\max}$. The $\mathcal{O}(1/\delta\sqrt{t})$ rate comes from the fact that both $\overline{L}_t$ and $h_0$ are upper bounded. $\overline{L}_t$ is bounded by Proposition~\ref{prop:lipschitz_estimate_bound} and $h_0$ is bounded by assumption.

\end{proof}

\subsection{Matching Pursuit}

In the context of Matching Pursuit, we propose the following criterion which we name the MP gap: $g^{\text{MP}}(\xx) = \max_{\ss \in \BB} \langle \nabla f(\xx), \ss\rangle$, where $\BB$ is as defined in \ref{apx:pseudocode_mp}.
Note that $g^{\text{MP}}$ is always non-negative and  $g^{\text{MP}}(\xx^{\star}) = 0$ implies $\langle \nabla f(\xx^{\star}), \ss\rangle = 0$ for all $\ss \in \BB$.
By linearity of the inner product we then have $\langle \nabla f(\xx^{\star}), \xx - \xx^\star \rangle = 0$ for any $\xx$ in the domain, since $\xx - \xx^\star$ lies in the linear span of $\mathcal{A}$.
Hence $\xx^{\star}$ is a stationary point and  $g^{\text{MP}}$ is an appropriate measure of stationarity for this problem.

\hfill

\hfill

\begin{mdframed}
  \begin{customtheorem}{\ref{thm:noncvx}.B}
    Let $\xx_{t}$ denote the iterate generated by \hyperref[alg:amp]{AdaMP} after $t$ iterations. Then for $t \geq 0$ we have the following suboptimality bound in terms of the MP gap:
    \begin{equation}
      \lim_{k\to \infty}g^{\text{MP}}(\xx_k) = 0\qquad\text{ and }\qquad \min_{0\leq k\leq t} g^{\text{MP}}(\xx_k)\leq \frac{\radius(\AA)}{\delta}\sqrt{ \frac{2h_0\overline{L}_t}{t+1}} = \mathcal{O}\left(\frac{1}{\delta \sqrt{t}}\right)~.
    \end{equation}
  \end{customtheorem}
\end{mdframed}

\begin{proof}
  The proof similar than that of Theorem~\ref{thm:noncvx}.A, except that in this case the expression of the step-size is simpler and does not depend on the minimum of two quantities. This avoids the case distinction that was necessary in the previous proof, resulting in a much simpler proof.

  For all $k= 0,\dots, t$, using the sufficient decrease condition, and the definitions of $\gamma_k$ and $g_k$:
  \begin{align}
    f(\xx_{k+1}) - f(\xx_k) & \leq \gamma_k \langle \nabla f(\xx_k), \dd_k\rangle + \frac{\gamma_k^2 L_{k}}{2}\|\dd_k\|^2                     \\
                            & \leq \min_{\eta \geq 0}\left\{-\eta g_k + \frac{1}{2}\eta^2 L_k\|\dd_k\|^2  \right\} \label{eq:mp_min_decrease} \\
                            & \leq -  \frac{g_k^2}{2L_k\|\dd_k\|^2}~,\label{eq:mp_decrease_solved}
  \end{align}
  where the last inequality comes from minimizing with respect to $\eta$.
  Summating over $k$ from $0$ to $t$ and negating the previous inequality, we obtain:
  \begin{align}
    \sum_{0\leq k\leq t}\frac{g_k^2}{L_k} & \leq (f(\xx_0) - f(\xx_t))\radius(\AA)^2\leq  2h_0\radius(\AA)^2~.
  \end{align}
  Taking the limit for $t \to \infty$ we obtain that the left hand side has a finite sum since the right hand side is bounded by assumption. Therefore, $g_k \to 0$, which by Lemma \ref{lemma:bound_gt_fw_gap} implies $\lim_{k\to \infty} g(\xx_k) = 0$. This proves the first claim of the Theorem.

  We now aim to derive explicit convergence rates. Taking the min over the $g_k$s and taking a square root for the last inequality
  \begin{align}
    \min_{0\leq k\leq t} g_k \leq \sqrt{ \frac{2h_0\radius(\AA)^2}{\sum_{0\leq k\leq t} {L_k}^{-1}}}
  \end{align}
  The term $\left({n}/{\sum_{0\leq k\leq t} {L_k}^{-1}}\right)$ is the \emph{harmonic mean} of the $L_k$s, which is always upper bounded by the average $\overline{L_t}$. Hence we obtain
  \begin{align}
    \min_{0\leq k\leq t} g_k \leq \frac{\radius(\AA)}{\delta}\sqrt{ \frac{2h_0 \overline{L}_t}{t+1}}~.
  \end{align}
  The claimed rate then follows from using the bound $g(\xx_k) \leq \frac{1}{\delta} g_k$ from Lemma~\ref{lemma:bound_gt_fw_gap}, valid for all $k\geq 0$.

  The $\mathcal{O}(1/\delta\sqrt{t})$ rate comes from the fact that both $\overline{L}_t$ and $h_0$ are upper bounded. $\overline{L}_t$ is bounded by Proposition~\ref{prop:lipschitz_estimate_bound} and $h_0$ is bounded by assumption.

\end{proof}

\paragraph{Note: Harmonic mean vs arithmetic mean.} The convergence rate for MP on non-convex objectives (Theorem~\ref{thm:noncvx}) also holds by replacing $\overline{L}_t$ by its harmonic mean $H_t \defas N_t/ ({\sum_{k=0}^{t-1}{L_k^{-1}}}\mathds{1}\{k \text{ is a good step}\})$ respectively. The harmonic mean is always less than the arithmetic mean, i.e., $H_t \leq \overline{L_t} $, although for simplicity we only stated both theorems with the arithmetic mean.
Note that the Harmonic mean is Schur-concave, implying that $H_t \leq t\min \{L_k: k\leq t\}$, i.e. it is controlled by the smallest Lipschitz estimate encountered so far.

\clearpage

\section{Proofs of convergence for convex objectives}\label{apx:proof_convex}

In this section we provide a proof the convergence rates stated in the theorem for convex objectives (Theorem \ref{thm:convex}). The section is structured as follows. We start by proving a technical result which is a slight variation of Lemma \ref{lemma:sufficient_decrease_inequality} and which will be used in the proof of Theorem \ref{thm:convex}. This is followed by the proof of Theorem~\ref{thm:convex}.

\subsection{Frank-Wolfe variants}

\begin{lemma}\label{lemma:decrease_good_step}
  For any of the proposed FW variants, if $t$ is a good step, then we have
  \begin{equation}
    f(\xx_{t+1}) \leq f(\xx_t)  -\xi {g_t}+ \frac{\xi^2 L_{t}}{2}\|\dd_t\|^2~ \text{ for all }\xi \in [0, 1].
  \end{equation}
\end{lemma}
\begin{proof}
  If $\gamma_t^{\max} \geq 1$, the result is obvious from Lemma \ref{lemma:sufficient_decrease_inequality}.
  If $\gamma_t^{\max} < 1$, then the inequality is only valid in the smaller interval $[0, \gamma_t^{\max}]$. However, since we have assumed that this is a good step, if $\gamma_t^{\max} < 1$ then we must have $\gamma_t < \gamma^{\max}_t$.
  By Lemma \ref{lemma:sufficient_decrease_inequality}, we have
  \begin{align}
    f(\xx_{t+1}) & \leq f(\xx_t) + \min_{\xi \in [0, \gamma_t^{\max}]}\left\{\xi\langle \nabla f(\xx_t), \dd_t\rangle + \frac{L_t\xi^2}{2}\|\dd_t\|^2 \right\}
  \end{align}
  Because $\gamma_t < \gamma_t^{\max}$ and since the expression inside the minimization term of the previous equation is a quadratic function of $\xi$, $\gamma_t$ is the unconstrained minimum and so we have
  \begin{align}
    f(\xx_{t+1}) & \leq f(\xx_t) + \min_{\xi \geq 0}\left\{\xi\langle \nabla f(\xx_t), \dd_t\rangle + \frac{L_t\xi^2}{2}\|\dd_t\|^2 \right\}       \\
                 & \leq f(\xx_t) + \min_{\xi \in [0, 1]}\left\{\xi\langle \nabla f(\xx_t), \dd_t\rangle + \frac{L_t\xi^2}{2}\|\dd_t\|^2 \right\}~.
  \end{align}
  The claimed bound then follows from the optimality of the min.
\end{proof}

The following lemma allows to relate the quantity $\langle \nabla f(\xx_t), \xx_t - \ss_t \rangle$ with a primal-dual gap and will be essential in the proof of Theorem~\ref{thm:convex}.
\begin{lemma}\label{lemma:primal_dual_relationship}
  Let $\ss_t$ be as defined in any of the FW variants. Then for any iterate $t \geq 0$ we have
  \begin{equation}
    \langle \nabla f(\xx_t), \xx_t - \ss_t \rangle \geq \delta(f(\xx_t) - \psi(\nabla f(\xx_t)))~.
  \end{equation}
\end{lemma}
\begin{proof}
  \begin{align}
    \langle \nabla f(\xx_t), \xx_t - \ss_t\rangle & \stackrel{\eqref{eq:lmo}}{\geq}  \delta\max_{\ss \in \conv(\mathcal{A})} \langle \nabla f(\xx_t), \xx_t - \ss\rangle                                                                        \\
                                                  & =  \delta \langle \nabla f(\xx_t), \xx_t\rangle + \delta \max_{\ss \in \conv(\mathcal{A})}\langle -\nabla f(\xx_t), \ss\rangle                                                              \\
                                                  & = \delta \big( \langle \nabla f(\xx_t), \xx_t\rangle + \sigma_{\conv(\mathcal{A})}(-\nabla f(\xx_t))\big)                                                                                   \\
                                                  & = \delta\big(f(\xx_t) + \underbrace{f^*(\nabla f(\xx_t)) + \sigma_{\conv(\mathcal{A})}(-\nabla f(\xx_t))}_{=-\psi(\nabla f(\xx_t))}\big) = \delta\big(f(\xx_t) - \psi(\nabla f(\xx_t))\big)
  \end{align}
  where the first identity uses the definition of $\ss_t$, the second one the definition of convex conjugate and the last one is a consequence of the Fenchel-Young identity. We recall $\sigma_{\conv(\mathcal{A})}$ is the support function of $\conv(\mathcal{A})$.

\end{proof}

\clearpage 

\begin{mdframed}
  \begin{customtheorem}{\ref{thm:convex}.A}
    Let $f$ be convex, $\xx_{t}$ denote the iterate generated by any of the proposed FW variants (\hyperref[alg:afw]{AdaFW}, \hyperref[alg:aafw]{AdaAFW}, \hyperref[alg:apfw]{AdaPFW}) after $t$ iterations, with $N_t \geq 1$, and let $\uu_t$ be defined recursively as $\uu_0 = \nabla f(\xx_0)$, $\uu_{t+1} = (1 - \xi_t)\uu_{t} + \xi_t\nabla f(\xx_t)$, where $\xi_t = 2 / (\delta N_t + 2)$ if $t$ is a good step and $\xi_t = 0$ otherwise. Then we have:
    \begin{equation}
      h_t \leq f(\xx_t) - \psi(\uu_t) \leq \frac{2\overline{L}_{t}\diam(\mathcal{A})^2}{ \delta^2 N_t + \delta} + \frac{2(1 - \delta )}{\delta^2 N_t^2 + \delta N_t} \big(f(\xx_0) - \psi(\uu_0)\big)= \mathcal{O}\left(\frac{1}{\delta^2 t}\right)~.
    \end{equation}
  \end{customtheorem}
\end{mdframed}
\begin{proof}
  The proof is structured as follows. First, we derive a bound for the case that $k$ is a good step. Second, we derive a bound for the case that $k$ is a bad step. Finally, we add over all iterates to derive the claimed bound.
  
  In both the good and bad step cases, we'll make use of the auxiliary variable $\sigma_k$. This is defined recursively as $\sigma_0 = \psi(\nabla f(\xx_k))$, $\sigma_{k+1} = (1 - \delta\xi_k)\sigma_k + \delta \xi_k \psi(\nabla f(\xx_k))$. Since $\psi$ is concave, Jensen's inequality gives that $\psi(\uu_k) \geq \sigma_k$ for all $k$. 

  \noindent \underline{\textbf{Case 1: $k$ is a good step}:}\\
  By  Lemma \ref{lemma:decrease_good_step}, we have the following sequence of inequalities, valid for all $\xi_t \in [0, 1]$:
  \begin{align}
    f(\xx_{k+1}) & \leq f(\xx_k) -\xi_k {g_k}+ \frac{\xi_k^2 L_{k}}{2}\|\dd_k\|^2                                            \\
                 & \leq  f(\xx_k) -\xi_k \langle \nabla f(\xx_k), \xx_k - \ss_k\rangle + \frac{\xi_k^2 L_{k}}{2}\|\dd_k\|^2  \\
                 & \leq (1 - \delta\xi_k)f(\xx_k) + \delta\xi_k \psi(\nabla f(\xx_k)) + \frac{\xi_k^2 L_{t}}{2}\|\dd_k\|^2~,
  \end{align}
  where the second inequality follows from the definition of $g_k$ (this is an equality for \hyperref[alg:afw]{AdaFW} but an inequality for the other variants) and the last inequality follows from Lemma \ref{lemma:primal_dual_relationship}.

  Subtracting $\psi(\uu_{k+1})$ from both sides of the previous inequality gives
  \begin{align}
  f(\xx_{k+1}) - \psi(\uu_{k+1}) &\leq f(\xx_{k+1}) - \sigma_{k+1} \\
    &\leq (1 - \delta\xi_k)\big[f(\xx_k) - \sigma_k\big] + \frac{\xi_k^2 L_k}{2}\|\ss_k - \xx_k\|^2\label{eq:pd_gap_tmp1}
  \end{align}
  Let $\xi_k = 2/ (\delta N_k + 2)$ and $a_k \defas \frac{1}{2}((N_k - 2)\delta + 2)((N_k - 1)\delta + 2)$. With these definitions, we have the following trivial identities that we will use soon:
  \begin{align}
    a_{k+1}(1 - \delta \xi_k) = \frac{1}{2}( (N_k - 2)\delta + 2)((N_k-1)\delta + 2) = a_k\label{eq:ak1} \\
    a_{k+1}\frac{\xi_k^2}{2} = \frac{((N_k-1)\delta + 2)}{(N_k\delta + 2)}\leq 1
  \end{align}
  where in the first inequality we have used that $k$ is a good step and so $N_{k+1} = N_k + 1$.

  Multiplying \eqref{eq:pd_gap_tmp1} by $a_{k+1}$ we have
  \begin{align}
    a_{k+1}\big(f(\xx_{k+1}) - \psi(\uu_{k+1})) & \leq a_{k+1}(1 - \delta\xi_k)\big[f(\xx_k) - \sigma_k\big] + \frac{L_k}{2}\|\ss_k - \xx_k\|^2       \\
                                             & \stackrel{\eqref{eq:ak1}}{=} a_k\big[f(\xx_k) - \sigma_k\big] + \frac{L_k}{2}\|\ss_k - \xx_k\|^2    \\
                                             & ~~{\leq }~\,\,a_k\big[f(\xx_k) - \sigma_k\big] + L_k \diam(\mathcal{A})^2\label{eq:bound_good_step}
  \end{align}

  \noindent \underline{\textbf{Case 2: $k$ is a bad step}:}\\
  Lemma~\ref{lemma:sufficient_decrease_inequality}  with $\xi_k=0$ guarantees that the objective function is non-increasing, i.e., $f(\xx_{k+1}) \leq f(\xx_k)$. By construction of $\sigma_k$ we have $\sigma_{k+1} = \sigma_k$, and so substracting both multiplied by $a_{k+1}$ we obtain
  \begin{align}
    a_{k+1}\big(f(\xx_{k+1}) - \psi(\uu_{k+1}) \big) & \leq a_{k+1}\big(f(\xx_{k+1}) - \sigma_{k+1} \big) \\
    & \leq a_{k+1}\big(f(\xx_{k}) - \sigma_{k} \big)                       \\
    & = a_{k}\big(f(\xx_{k}) - \sigma_{k} \big)~,\label{eq:bound_bad_step}
  \end{align}
  where in the last identity we have used that its a bad step and so $a_{k+1} = a_k$.

  \noindent \underline{\textbf{Final: combining cases and adding over iterates}:}\\
  We can combine \eqref{eq:bound_good_step} and \eqref{eq:bound_bad_step} into the following inequality:
  \begin{equation}
    a_{k+1}\big(f(\xx_{k+1}) - \psi(\uu_{k+1}) \big) - a_k \big(f(\xx_{k}) - \psi(\uu_{k}) \big) \leq L_{k}\diam(\AA)^2\mathds{1}\{k \text{ is a good step}\}~,
  \end{equation}
  where $\mathds{1}\{\text{condition}\}$ is 1 if condition is verified and 0 otherwise.

  Adding this inequality from $0$ to $t - 1$ gives
  \begin{align}
    a_t\big(f(\xx_{t}) - \psi(\uu_{t}) \big) & \leq \sum_{k = 0}^{t - 1} L_{k}Q_{A}^2\mathds{1}\{k \text{ is a good step}\} + a_0(f(\xx_{0}) - \sigma_0) \\
                                          & = N_{t} \overline{L}_{t}\diam(\mathcal{A})^2 + {(1 - \delta )(2 - \delta)} (f(\xx_{0}) - \sigma_{0})
  \end{align}
  Finally, dividing both sides by $a_{t}$ (note that $a_t > 0$ for $N_t \geq 1$) and using $(2 - \delta) \leq 2$ we obtain
  \begin{align}
    f(\xx_{t}) - \psi(\uu_{t}) & \leq \frac{2 N_t}{ ((N_t-2)\delta + 2)((N_t - 1)\delta + 2)}  \overline{L}_{t}Q_{\mathcal{A}}^2  \\
                          & \qquad + \frac{4(1 - \delta )}{ ((N_t-2)\delta + 2)((N_t - 1)\delta + 2)}  (f(\xx_0) - \sigma_0)
  \end{align}
  We will now use the inequalities $(N_t-2)\delta + 2 \geq N_t\delta$ and $(N_t - 1)\delta + 2 \geq N_t\delta + 1$ for the terms in the denominator to obtain
  \begin{align}
     f(\xx_{t}) - \psi(\uu_{t}) & \leq \frac{2 \overline{L}_{t}Q_{\mathcal{A}}^2}{ \delta^2 N_t + \delta}   + \frac{4(1 - \delta )}{\delta_t^2 N_t^2 + \delta N_t}  (f(\xx_0) - f(\xx^\star))~.\label{eq:sublinear_bound_sigma}
  \end{align}
 which is the desired bound:
  \begin{equation}
    f(\xx_{t}) - \psi(\uu_{t}) \leq \frac{2 \overline{L}_{t}Q_{\mathcal{A}}^2}{ \delta^2 N_t + \delta}   + \frac{4(1 - \delta )}{\delta_t^2 N_t^2 + \delta N_t}  \big[f(\xx_0) - \psi(\nabla f(\xx_0))\big].
  \end{equation}

  We will now prove the bound $h_t \leq f(\xx_{t}) - \psi(\uu_{t})$. Let $\uu^\star$ be an arbitrary maximizer of $\psi$. Then by duality we have that $f(\xx^\star) = \psi(\uu^\star)$ and so
  \begin{align}
    f(\xx_{t}) - \psi(\uu_{t}) = f(\xx_{t}) - f*(\xx^\star) + \psi(\uu^\star) - \psi(\uu_{t}) \geq f(\xx_{t}) - f*(\xx^\star) = h_t
  \end{align}

  Finally, the $\mathcal{O}(\frac{1}{\delta t})$ rate comes from bounding the number of good steps from \eqref{eq:bad_steps}, for which we have $1/N_t \leq \mathcal{O}(1/t)$, and bounding the Lipschitz estimate by a contant (Proposition~\ref{prop:lipschitz_estimate_bound}).
\end{proof}

\subsection{Matching Pursuit}

\begin{lemma}\label{lemma:ineq_mp_convex}
  Let $\ss_t$ be as defined in \hyperref[alg:amp]{AdaMP}, $R_\BB$ be the level set radius defined as
  \begin{equation}
    R_\BB = \max_{\substack{\xx \in \lin(\AA)\\f(\xx) \leq f(\xx_0)}} \|\xx - \xx^\star\|_{\BB}~,
  \end{equation}
  and $\xx^\star$ be any solution to \eqref{eq:obj_fun_mp}. Then we have
  \begin{equation}
    \langle -\nabla f(\xx_t), \ss_t\rangle \geq \frac{\delta}{\max\{\mathcal{R}_\BB, 1\}}\big(f(\xx_t) - f(\xx^\star)\big)
  \end{equation}
\end{lemma}
\begin{proof}
  By definition of atomic norm we have
  \begin{equation}
    \frac{\xx_t - \xx^\star_t}{\|\xx_t - \xx^\star\|_\BB} \in \conv(\BB)
  \end{equation}
  Since $f(\xx_t) \leq f(\xx_0)$, which is a consequence of sufficient decrease condition (Eq.~\eqref{eq:mp_min_decrease}), we have that $R_\BB\geq \|\xx_t - \xx^\star\|_\BB$ and so $\zeta \defas \|\xx_t - \xx^\star\|_\BB/R_\BB \leq 1$. By symmetry of $\BB$ we have that
  \begin{equation}\label{eq:diff_in_B}
    \frac{\xx_t - \xx^\star}{R_\BB} = \zeta
    \frac{\xx_t - \xx^\star}{\|\xx_t - \xx^\star\|_\BB} + (1 - \zeta)\boldsymbol{0} \in \conv(\BB)~.
  \end{equation}

  We will now use this fact to bound the original expression. By definition of $\ss_t$ we have
  \begin{align}
    \langle -\nabla f(\xx_t), \ss_t\rangle & \stackrel{\eqref{eq:mp_lmo}}{\geq} \delta\max_{\ss \in \BB}\langle -\nabla f(\xx_t), \ss\rangle~             \\
                                           & \stackrel{\eqref{eq:diff_in_B}}{\geq}\frac{\delta}{R_\BB}\langle -\nabla f(\xx_t), \xx_t - \xx^\star\rangle~ \\
                                           & ~\geq \frac{\delta}{R_\BB}(f(\xx_t) - f(\xx^\star))
  \end{align}
  where the last inequality follows by convexity.

\end{proof}

\begin{mdframed}
  \begin{customtheorem}{\ref{thm:convex}.B}
    Let $f$ be convex, $\xx^\star$ be an arbitrary solution to \eqref{eq:obj_fun_mp} and let $R_\BB$ the level set radius:
    \begin{equation}
      R_\BB = \max_{\substack{\xx \in \lin(\AA)\\f(\xx) \leq f(\xx_0)}} \|\xx - \xx^\star\|_{\BB}~.
    \end{equation}
    If we denote by $\xx_{t}$ the iterate generated by \hyperref[alg:amp]{AdaMP} after $t \geq 1$ iterations and $\beta = \delta / R_\BB$, then we have:
    \begin{equation}
      f(\xx_t) - f(\xx^\star) \leq \frac{2\overline{L}_{t}\radius(\AA)^2}{ \beta^2 t + \beta} + \frac{2(1 - \beta )}{\beta^2 t^2 + \beta t} h_0= \mathcal{O}\left(\frac{1}{\beta^2 t}\right)~.
    \end{equation}
  \end{customtheorem}
\end{mdframed}
\begin{proof}
  Let $\xx^\star$ be an arbitrary solution to \eqref{eq:obj_fun_mp}. Then by Lemma \ref{lemma:sufficient_decrease_inequality}, we have the following sequence of inequalities, valid for all $\xi_t \geq 0$:
  \begin{align}
    f(\xx_{k+1}) & \leq f(\xx_k)  -\xi_k \langle  -\nabla f(\xx_k), \ss_k\rangle + \frac{\xi_k^2 L_{k}}{2}\|\ss_k\|^2                 \\
                 & \leq  f(\xx_k) - \xi_k\frac{\delta}{R_\BB}\big[f(\xx_k) - f(\xx^\star)\big] + \frac{\xi_k^2 L_{t}}{2}\|\ss_k\|^2~,
  \end{align}
  where the second inequality follows from Lemma \ref{lemma:ineq_mp_convex}.

  Subtracting $f(\xx^\star)$ from both sides of the previous inequality gives
  \begin{align}\label{eq:pd_gap_tmp1_mp}
    f(\xx_{k+1}) - f(\xx^\star) \leq \left(1 - \frac{\delta}{R_\BB}\xi_k\right)\big[f(\xx_k) - f(\xx^\star)\big] + \frac{\xi_k^2 L_k}{2}\|\ss_k\|^2~.
  \end{align}

  Let $\beta={\delta}/{R_\BB}$ and $\xi_k = 2/ (\beta k + 2)$ and $a_k \defas \frac{1}{2}((k - 2)\beta + 2)((k - 1)\beta + 2)$. With these definitions, we have the following trivial results:
  \begin{align}
    a_{k+1}(1 - \beta \xi_k) = \frac{1}{2}( (k - 2)\beta + 2)((k-1)\beta + 2) = a_k\label{eq:ak1_mp} \\
    a_{k+1}\frac{\xi_k^2}{2} = \frac{((k-1)\beta + 2)}{(k\beta + 2)}\leq 1~.
  \end{align}

  Multiplying \eqref{eq:pd_gap_tmp1_mp} by $a_{k+1}$ we have
  \begin{align}
    a_{k+1}\big(f(\xx_{k+1}) - f(\xx^\star)) & \leq a_{k+1}(1 - \beta\xi_k)\big[f(\xx_k) - f(\xx^\star)\big] + \frac{L_k}{2}\|\ss_k\|^2                        \\
                                             & \stackrel{\eqref{eq:ak1}}{=} a_k\big[f(\xx_k) - f(\xx^\star)\big] + \frac{L_k}{2}\|\ss_k\|^2                    \\
                                             & \quad{\leq }~\,\,a_k\big[f(\xx_k) - f(\xx^\star)\big] + L_t \radius(\mathcal{A})^2\label{eq:bound_good_step_mp}
  \end{align}

  Adding this last inequality from $0$ to $t - 1$ gives
  \begin{align}
    a_t\big(f(\xx_{t}) - f(\xx^\star) \big) & \leq \sum_{k=0}^{t-1} L_{k}\radius(\AA)^2 + a_0(f(\xx_{0}) - \beta_0)                           \\
                                            & = t \overline{L}_{t}\diam(\mathcal{A})^2 + {(1 - \delta )(2 - \delta)} (f(\xx_{0}) - \beta_{0})
  \end{align}
  Finally, dividing both sides by $a_{t}$ (note that $a_1 = 2 -\beta \geq 1$ and so $a_t$ is strictly positive for $t \geq 1$),  and using $(2 - \delta) \leq 2$ we obtain
  \begin{align}
    f(\xx_{t}) - f(\xx^\star) & \leq \frac{2 t}{ ((t-2)\beta + 2)((t - 1)\beta + 2)}  \overline{L}_{t}\radius(\AA)^2     \\
                              & \qquad + \frac{4(1 - \beta )}{ ((t-2)\beta + 2)((t - 1)\beta + 2)}  (f(\xx_0) - \beta_0)
  \end{align}
  We will now use the inequalities $(t-2)\beta + 2 \geq t\beta$ and $(t - 1)\beta + 2 \geq t\beta + 1$ to simplify the terms in the denominator. With this we obtain to obtain
  \begin{align}
    f(\xx_{t}) - f(\xx^\star) & \leq \frac{2 \overline{L}_{t}\radius(\AA)^2}{ \beta^2 N_t + \beta}   + \frac{4(1 - \beta )}{\beta_t^2 N_t^2 + \beta N_t}  (f(\xx_0) - f(\xx^\star))~,\label{eq:sublinear_bound_sigma_mp}
  \end{align}
  which is the desired bound.
\end{proof}

\clearpage

\section{Proofs of convergence for strongly convex objectives}\label{apx:linear_proofs}

The following proofs depend on some definitions of geometric constants, which are defined in \ref{apx:basic_defs} as well as two crucial lemmas from \citep[Appendix C]{lacoste2015global}.

\subsection{Frank-Wolfe variants}
We are now ready to present the convergence rate of the backtracking Frank--Wolfe variants. As we did in \ref{apx:proof_nonconvex}, although the original proof combines the rates for FW variants and MP, the proof will be split into two, in which we prove separately the linear convergence rates for \hyperref[alg:aafw]{AdaAFW} and \hyperref[alg:apfw]{AdaPFW} (Theorem~\ref{thm:linear_conv_afw}.A) and \hyperref[alg:amp]{AdaMP} (Theorem~\ref{thm:linear_conv_afw}.B).

\begin{mdframed}
  \begin{customtheorem}{\ref{thm:linear_conv_afw}.A}
    Let $f$ be $\mu$--strongly convex. Then for each good step we have the following geometric decrease:
    \begin{equation}
      h_{t+1} \leq (1 - \rho_t) h_t,
    \end{equation}
    with
    \begin{align}
      \rho_t & = \dfrac{\mu \delta^2}{4L_t}\left(\dfrac{\PWidth(\AA)}{\diam(\conv(\mathcal{A}))}\right)^2                                      & \text{for \hyperref[alg:aafw]{AdaAFW}} \\
      \rho_t & = \min\Big\{ \dfrac{\delta}{2}, \delta^2\dfrac{\mu }{L_t} \left(\dfrac{\PWidth(\AA)}{\diam(\conv(\mathcal{A}))}\right)^2 \Big\} & \text{for \hyperref[alg:apfw]{AdaPFW}}
    \end{align}
  \end{customtheorem}
\end{mdframed}

\paragraph{Note.} In the main paper we provided the simplified bound $\rho_t = \dfrac{\mu}{4{L}_t}\left(\dfrac{\PWidth(\AA)}{\diam(\AA)}\right)^2$ for both algorithms \hyperref[alg:aafw]{AdaAFW} and \hyperref[alg:apfw]{AdaPFW} for simplicity. It is easy to see that the bound for \hyperref[alg:apfw]{AdaPFW} above can be trivially bounded by this quantity by noting that $\delta^2 \leq \delta$ and that $\mu/L_t$ and ${\PWidth(\AA)}/{\diam(\conv(\mathcal{A}))}$ are necessarily smaller than 1.

\begin{proof}
  The structure of this proof is similar to that of \citep[Theorem 8]{lacoste2015global}. We begin by upper bounding the suboptimality $h_t$. Then we derive a lower bound on $h_{t+1} - h_t$. Combining both we arrive at the desired geometric decrease.\\

  \noindent \underline{\emph{Upper bounding $h_t$}}\\
  Assume $\xx_t$ is not optimal, ie $h_t > 0$. Then we have $\langle -\nabla f(\xx_t), \xx^\star - \xx_t \rangle > 0$. Using the definition of the geometric strong convexity bound and letting $\overline{\gamma} \defas \gamma (\xx_t, \xx^\star)$ we have
  \begin{align}
    \dfrac{\overline{\gamma}^2}{2} \mu_f^A & \leq f(\xx^\star) - f(\xx_t) + \langle - \nabla f(\xx_t), \xx^\star - \xx_t \rangle       \\
                                           & = -h_t + \overline{\gamma} \langle - \nabla f(\xx_t), \ss_f(\xx_t) - \vv_f(\xx_t) \rangle \\
                                           & \leq -h_t + \overline{\gamma} \langle - \nabla f(\xx_t), \ss_t - \vv_t  \rangle           \\
                                           & = -h_t + \overline{\gamma} q_t~,
  \end{align}
  where $q_t \defas \langle -\nabla f(\xx_t), \ss_t - \vv_t\rangle$. For the last inequality we have used the definition of $\vv_f(\xx)$ which implies $\langle f(\xx_t), \vv_f(\xx_t) \rangle \leq \langle \nabla f(\xx_t), \vv_t \rangle$ and the fact that $\ss_t = \ss_f(\xx_t)$. Therefore
  \begin{align}
    h_t \leq -\dfrac{\overline{\gamma}^2}{2} \mu_f^A + \overline{\gamma} q_t~, \label{ineq:upperbound}
  \end{align}
  which can always be upper bounded by taking $\overline{\gamma} = \mu^{-1} q_t$ (since this value of $\overline{\gamma}$ maximizes the expression on the right hand side of the previous inequality) to arrive at
  \begin{align}
    h_t & \leq \frac{q_t^2}{2\mu_f^A}                                   \\
        & \leq \frac{q_t^2}{2\mu \Delta^2}~,\label{eq:strong_convexity}
  \end{align}
  with $\Delta \defas \PWidth(\AA)$ and where the last inequality follows from Lemma \ref{lemma:lb_mufa}.

  \underline{\emph{Lower bounding progress $h_t - h_{t+1}$.}}\\
  Let $G$ be defined as $G = 1/2$ for \hyperref[alg:aafw]{AdaAFW} and $G = 1$ for \hyperref[alg:apfw]{AdaPFW}.
  We will now prove that for both algorithms we have
  \begin{equation}\label{eq:gap_inequality}
    \langle -\nabla f(\xx_t), \dd_t\rangle \geq
    \delta G q_t~.
  \end{equation}
  For \hyperref[alg:aafw]{AdaAFW}, by the way the direction $\dd_t$ is chosen on Line \ref{line:decision_step}, we have the following sequence of inequalities:
  \begin{align*}
    2 \langle -\nabla f(\xx_t) , \dd_t \rangle
     & \geq \langle -\nabla f(\xx_t) , \dd_t^{FW} \rangle + \langle -\nabla f(\xx_t) , \dd_t^{A} \rangle                       \\
     & \geq \delta \langle  -\nabla f(\xx_t) , \ss_t - \xx_t \rangle + \delta \langle -\nabla f(\xx_t) , \xx_t - \vv_t \rangle \\
     & = \delta \langle -\nabla f(\xx_t) , \ss_t - \vv_t \rangle                                                               \\
     & = \delta q_t~,
  \end{align*}
  For \hyperref[alg:apfw]{AdaPFW}, since $\dd_t = \ss_t - \vv_t$, it follows from the definition of $q_t$ that $\langle - \nabla f(\xx_t), \dd_t \rangle \geq \delta q_t$.

  We split the rest of the analysis into three cases: $\gamma_t < \gamma_t^{\max}, \gamma_t = \gamma_t^{\max} \geq 1$ and $ \gamma_t = \gamma_t^{\max} < 1$. We prove a geometric descent in the first two cases. In the case where $\gamma_t = \gamma_t^{\max} < 1$ (a bad step) we show that the number of bad steps is bounded. \\

  \noindent \underline{\textbf{Case 1: $\gamma_t < \gamma_t^{\max}$}:}\\
  By Lemma \ref{lemma:sufficient_decrease_inequality}, we have
  \begin{align}
    f(\xx_{t+1}) = f(\xx_t + \gamma_t \dd_t) & \leq f(\xx_t) + \min_{\eta \in [0, \gamma_t^{\max}]}\left\{\eta\langle \nabla f(\xx_t), \dd_t\rangle + \frac{L_t\eta^2}{2}\|\dd_t\|^2 \right\} \label{eq:min-suffdec}
  \end{align}
  Because $\gamma_t < \gamma_t^{\max}$ and since the expression inside the minimization term \eqref{eq:min-suffdec} is a convex function of $\eta$, the minimizer is unique and it coincides with the minimum of the unconstrained problem. Hence we have
  \begin{align}\label{eq:dom-eq}
    \min_{\eta \in [0, \gamma_t^{\max}]}\left\{\eta\langle \nabla f(\xx_t), \dd_t\rangle + \frac{L_t\eta^2}{2}\|\dd_t\|^2 \right\} = \min_{\eta \geq 0}\left\{\eta\langle \nabla f(\xx_t), \dd_t\rangle + \frac{L_t\eta^2}{2}\|\dd_t\|^2 \right\}
  \end{align}
  Replacing in \eqref{lemma:sufficient_decrease_inequality}, our bound becomes
  \begin{align}
    f(\xx_{t+1}) = f(\xx_t + \gamma_t \dd_t) & \leq f(\xx_t) + \min_{\eta \geq 0}\left\{\eta\langle \nabla f(\xx_t), \dd_t\rangle + \frac{L_t\eta^2}{2}\|\dd_t\|^2 \right\} \\
                                             & \leq f(\xx_t) + \min_{\eta \geq 0}\left\{\eta\langle \nabla f(\xx_t), \dd_t\rangle + \frac{L_t\eta^2}{2} M^2 \right\}        \\
                                             & \leq f(\xx_t) + \eta\langle \nabla f(\xx_t), \dd_t\rangle + \frac{L_t\eta^2}{2} M^2 , \;\; \forall \eta \geq 0
  \end{align}
  where the second inequality comes from bounding $\|\dd_t\|$ by $M \defas \diam (\conv(\mathcal{A}))$. Subtracting $f(\xx^\star)$ from both sides and rearranging we have
  \begin{align}
    h_t - h_{t+1} & \geq \eta \langle -\nabla f(\xx_t), \dd_t \rangle - \dfrac{1}{2} \eta^2 L_t M^2, \;\; \forall \eta \geq 0~.
  \end{align}
  Using the gap inequality \eqref{eq:gap_inequality} our lower bound becomes
  \begin{align}\label{eq:lower-bound-gamma}
    h_t - h_{t+1} & \geq \eta \delta G q_t - \dfrac{1}{2} \eta^2 L_t M^2, \;\; \forall \eta \geq 0~.
  \end{align}

  Noting that the lower bound in \eqref{eq:lower-bound-gamma} is a concave function of $\eta$, we maximize the bound by selecting $\eta^\star = (L_t M^2)^{-1} \delta G q_t$. Plugging $\eta^\star$ into the bound in \eqref{eq:lower-bound-gamma} and then using the strong convexity bound \eqref{eq:strong_convexity} we have
  \begin{align}
    h_t - h_{t+1} \geq \dfrac{\mu G^2 \Delta^2 \delta^2}{L_t M^2} h_t \Longrightarrow h_{t+1} \leq \Big(1 - \dfrac{\mu G^2 \Delta^2 \delta^2}{L_t M^2}\Big) h_t~. \label{eq:lower-bound-strongcvx}
  \end{align}
  Then we have geometric convergence with rate $1-\rho$ where $\rho = (4L_tM^2)^{-1} \mu \Delta^2 \delta^2$ for \hyperref[alg:aafw]{AdaAFW} and $\rho = (L_t M^2)^{-1} \mu \Delta^2 \delta^2$ for \hyperref[alg:apfw]{AdaPFW}.\\

  \noindent \underline{\textbf{Case 2: $\gamma_t = \gamma_t^{\max} \geq 1$}} \\
  By Lemma~\ref{lemma:sufficient_decrease_inequality} and the gap inequality \eqref{eq:gap_inequality}, we have
  \begin{align}
    h_t - h_{t+1} = f(\xx_t) - f(\xx_{t+1})\geq \eta \delta G q_t - \dfrac{1}{2} \eta^2 L_t M^2, \; \; \forall \eta \leq \gamma_t^{\max}~. \label{eq:lower-bound-gamma-case2}
  \end{align}
  Since the lower bound \eqref{eq:lower-bound-gamma-case2} is true for all $\eta \leq \gamma_t^{\max}$, we can maximize the bound with $\eta^\star = \min \{ (L_tM^2)^{-1} \delta G q_t, \gamma_t^{\max}\}$. In the case when $\eta^\star = (L_tM^2)^{-1} \delta G q_t$ we get the same bound as we do in \eqref{eq:lower-bound-strongcvx} and hence have linear convergence with rate $1 - \rho$ where $\rho = (4L_tM^2)^{-1} \mu \Delta^2 \delta^2$ for \hyperref[alg:aafw]{AdaAFW} and $\rho = (L_t M^2)^{-1} \mu \Delta^2 \delta^2$ for \hyperref[alg:apfw]{AdaPFW}. If $\eta^\star = \gamma_t^{\max}$ then this implies $L_t M^2 \leq \delta G q_t$. Since $\gamma_t^{\max}$ is assumed to be greater than $1$ and the bound holds for all $\eta \leq \gamma_t^{\max}$ we have in particular that it holds for $\eta = 1$ and hence
  \begin{align}
    h_t - h_{t+1} & \geq \delta G q_t - \dfrac{1}{2} L_t M^2    \\
                  & \geq \delta G q_t - \dfrac{\delta G q_t}{2} \\
                  & \geq \dfrac{\delta G h_t}{2}~,
  \end{align}
  where in the second line we use the inequality $L_t M^2 \leq \delta G q_t$ and in the third we use the inequality $h_t \leq q_t$ which is an immediate consequence of convexity of $f$. Then we have
  \begin{align}
    h_{t+1} \leq (1 - \rho) h_t~,
  \end{align}
  where $\rho = \delta/4$ for \hyperref[alg:aafw]{AdaAFW} and $\rho = \delta/2$ for \hyperref[alg:apfw]{AdaPFW}. Note by Proposition \ref{prop:ub:pwdith} and the fact $\mu \leq L_t$ we have $\delta/4 \geq (4L_t M^2)^{-1} \mu \Delta^2 \delta^2$.\\

  \noindent \underline{\textbf{Case 3: $\gamma_t = \gamma_t^{\max} < 1$ (bad step)}}\\
  In this case, we have either a drop or swap step and can make no guarantee on the progress of the algorithm (drop and swap are defined in \ref{apx:sufficient_decrease}). For \hyperref[alg:aafw]{AdaAFW}, $\gamma_t = \gamma_t^{\max} < 1$ is a drop step. From lines \ref{line:decision_step}--\ref{l:gammamax2} of \hyperref[alg:aafw]{AdaAFW} we can make the following distinction of cases. In case of a  FW step, then $\mathcal{S}_{t+1} = \{ \ss_t\}$ and $\gamma_t = \gamma_t^{\max} = 1$, otherwise $\mathcal{S}_{t+1} = \mathcal{S}_{t} \cup \{\ss_t\}$. In case of an Away step, $\mathcal{S}_{t+1} = \mathcal{S}_{t} \backslash \{ \vv_t\}$ if $\gamma_t = \gamma_t^{\max} < 1$ , otherwise $\mathcal{S}_{t+1} = \mathcal{S}_{t}$. Note a drop step can only occur at an Away step. For \hyperref[alg:apfw]{AdaPFW}, $\gamma_t = \gamma_t^{\max} < 1$ will be a drop step when $\ss_t \in \mathcal{S}_{t}$ and will be a swap step when $\ss_t \not \in \mathcal{S}_{t}$.

  Even though at these bad steps we do not have the same geometric decrease, Lemma~\ref{lemma:drop_steps} yields that the sequence $\{h_t\}$ is a non-increasing sequence, i.e., $h_{t+1} \leq h_t$. Since we are guaranteed a geometric decrease on steps that are not bad steps, the bounds on the number of bad steps of Eq.~\eqref{eq:bad_steps} is sufficient to conclude that \hyperref[alg:aafw]{AdaAFW} and \hyperref[alg:apfw]{AdaPFW} exhibit a global linear convergence.

\end{proof}

\subsection{Matching Pursuit}

We start by proving the following lemma, which will be crucial in the proof of the backtracking MP's linear convergence rate.

\begin{lemma}\label{lemma:strong_cvx_lowerbound}
  Suppose that $\AA$ is a non-empty compact set and that $f$ is $\mu$--strongly convex.
  Let $\nabla_\BB f(\xx)$ denote the orthogonal projection of $\,\nabla f(\xx)$ onto $\lin(\BB)$. Then for all $\xx^\star - \xx \in \lin(\AA)$, we have
  \begin{equation}
    f(\xx^\star) \geq f(\xx) - \dfrac{1}{2\mu \mDW(\BB)^2} \| \nabla_{\BB} f(\xx)\|_{\BB^\star}^2~.
  \end{equation}
\end{lemma}

\begin{proof}


  From \citet[Theorem 6]{locatello2018revisiting}, we have that if $f$ is $\mu$-strongly convex, then
  \begin{align}
    \mu_{\BB} \defas \inf_{\xx,\yy \in \lin(\BB), \xx \not= \yy} \dfrac{2}{\|\yy - \xx\|_{\BB}^2} [f(\yy) - f(\xx) - \langle \nabla f(\xx), \yy - \xx \rangle]
  \end{align}
  is positive and verifies  $\mu_\BB \geq \mDW(\BB)^2 \mu$. Replacing $\yy = \xx + \gamma(\xx^\star - \xx)$ in the definition above we have
  \begin{equation}
    f(\xx + \gamma (\xx^\star - \xx)) \geq f(\xx) + \gamma \langle \nabla f(\xx), \xx^\star - \xx \rangle + \gamma^2 \dfrac{\mu_\BB}{2} \|\xx^\star - \xx\|_\BB^2 ~.
  \end{equation}
  We can fix $\gamma = 1$ on the left hand side and since the expression on the right hand side is true for all $\gamma$, we minimize over $\gamma$ to find $\gamma^* = -\langle \nabla f(\xx), \xx^\star - \xx \rangle/ \mu_\BB \|\xx^\star - \xx\|_\BB^2$. Thus the lower bound becomes
  \begin{align}
    f(\xx^\star) & \geq f(\xx) - \dfrac{1}{2\mu_\BB} \dfrac{\langle \nabla f(\xx), \xx^\star - \xx \rangle}{\| \xx^\star - \xx\|_\BB^2}            \\
                 & \geq f(\xx) - \dfrac{1}{2\mu \mDW(\BB)^2} \dfrac{\langle \nabla f(\xx), \xx^\star - \xx \rangle}{\| \xx^\star - \xx\|_{\BB}^2}  \\
                 & = f(\xx) - \dfrac{1}{2\mu \mDW(\BB)^2} \dfrac{\langle \nabla_\BB f(\xx), \xx^\star - \xx \rangle}{\| \xx^\star - \xx\|_{\BB}^2} \\
                 & \geq f(\xx) - \dfrac{1}{2\mu \mDW(\BB)^2} \| \nabla_\BB f(\xx) \|_{\BB^*}^2~,
  \end{align}
  where the last inequality follows by $|\langle \yy, \zz \rangle| \leq \| \yy\|_{\BB^*} \|\zz\|_\BB$
\end{proof}

\begin{mdframed}
  \begin{customtheorem}{\ref{thm:linear_conv_afw}.B} (Convergence rate backtracking MP)
    Let $f$ be $\mu$--strongly convex and suppose $\BB$ is a non-empty compact set. Then \hyperref[alg:amp]{AdaMP} verifies the following geometric decrease for each $t \geq 0$:
    \begin{equation}
      h_{t+1} \leq \Big( 1 - \delta^2\rho_t\Big) h_t, \;\; \text{ with } \rho_t = \dfrac{\mu}{{L}_t} \left(\dfrac{\mDW(\BB)}{\radius(\BB)}\right)^2~,
    \end{equation}
    where $\mDW(\BB)$ the minimal directional width of $\BB$.
  \end{customtheorem}
\end{mdframed}

\begin{proof}
  By Lemma~\ref{lemma:sufficient_decrease_inequality} and bounding $\|\dd_t\|$ by $R = \radius (\BB)$ we have
  \begin{align}
    f(\xx_{t+1}) & \leq f(\xx_t) + \min_{\eta \in \RR} \Big\{ \eta \langle \nabla f(\xx_t), \ss_t\rangle + \dfrac{\eta^2 L_t R^2}{2} \Big\} \\
                 & = f(\xx_t)  - \dfrac{\langle \nabla f(\xx_t), \ss_t \rangle^2}{2L_t R^2}                                                 \\
                 & \leq f(\xx_t) - \delta^2 \dfrac{\langle \nabla f(\xx_t), \ss_t^\star \rangle^2}{2L_t R^2} \label{eq:parallel_decomp}
  \end{align}

  where $\ss_t^\star$ is any element such that $\ss_t^\star \in \argmin_{\ss \in \BB} \langle \nabla f(\xx_t), \ss \rangle$ and the inequality follows from the optimality of min and the fact that $\langle \nabla f(\xx_t), \ss_t^\star \rangle \leq 0$. Let $\nabla_\BB f(\xx_t)$ denote as in Lemma~\ref{lemma:strong_cvx_lowerbound} the orthogonal projection of $\nabla f(\xx_t)$ onto $\lin(\BB)$. Then the previous inequality simplifies to
  \begin{align}
    f(\xx_{t+1}) \leq f(\xx_t) - \delta^2 \dfrac{\langle \nabla_\BB f(\xx_t), \ss_t^\star \rangle^2}{2L_t R^2}~. \label{eq:parallel_dual}
  \end{align}
  By definition of dual norm, we also have $\langle -\nabla_\BB f(\xx_t), \ss_t^\star \rangle = \| \nabla_\BB f(\xx_t) \|_{\BB*}^2$. Subtracting $f(\xx^\star)$ from both sides we obtain the upper-bound:
  \begin{align}
    h_{t+1} \leq h_t - \delta^2 \dfrac{\| \nabla_\BB f(\xx_t)\|_{\BB*}^2}{2L_t R^2}
  \end{align}

  To derive the lower-bound, we use Lemma \ref{lemma:strong_cvx_lowerbound} with $\xx = \xx_t$ and see that
  \begin{align}
    \| \nabla_\BB f(\xx_t) \|_{\BB^*} \geq 2 \mu \mDW(\BB)^2 h_t
  \end{align}

  Combining the upper and lower bound together we have
  \begin{align}
    h_{t+1} \leq \Big( 1 - \delta^2 \dfrac{\mu \mDW(\BB)^2}{L_t  R^2} \Big) h_t ~,
  \end{align}
  which is the claimed bound.
\end{proof}

\clearpage




\section{Experiments}

In this appendix we give give some details on the experiments which were omitted from the main text, as well as an extended set of results.

\subsection{$\ell_1$-regularized logistic regression, Madelon dataset}

For the first experiment, we consider an $\ell_1$-regularized logistic regression of the form
\begin{equation}\label{eq:l1_logistic_regression_apx}
  \argmin_{\|\xx\|_1 \leq \beta} \ \frac{1}{n} \sum_{i=1}^n \varphi(\aa_i^\top \xx, \boldsymbol{b}_i) + \frac{\lambda}{2} \|\xx\|^2_2~,
\end{equation}
where $\varphi$ is the logistic loss. The linear subproblems in this case can be computed exactly ($\delta=1$) and consists of finding the largest entry of the gradient. The regularization parameter $\lambda$ is always set to $\lambda=\frac{1}{n}$.

We first consider the case in which the data $\aa_i, \boldsymbol{b}_i$ is the Madelon datset.
Below are the curves objective suboptimality vs time for the different methods considered. The regularization parmeter, denoted $\ell_1$ ball radius in the figure, is chosen as to give $1\%$, $5\%$ and $20\%$ of non-zero coefficients (the middle figure is absent from the main text).

\begin{figure}[ht]
  \centering\includegraphics[width=\linewidth]{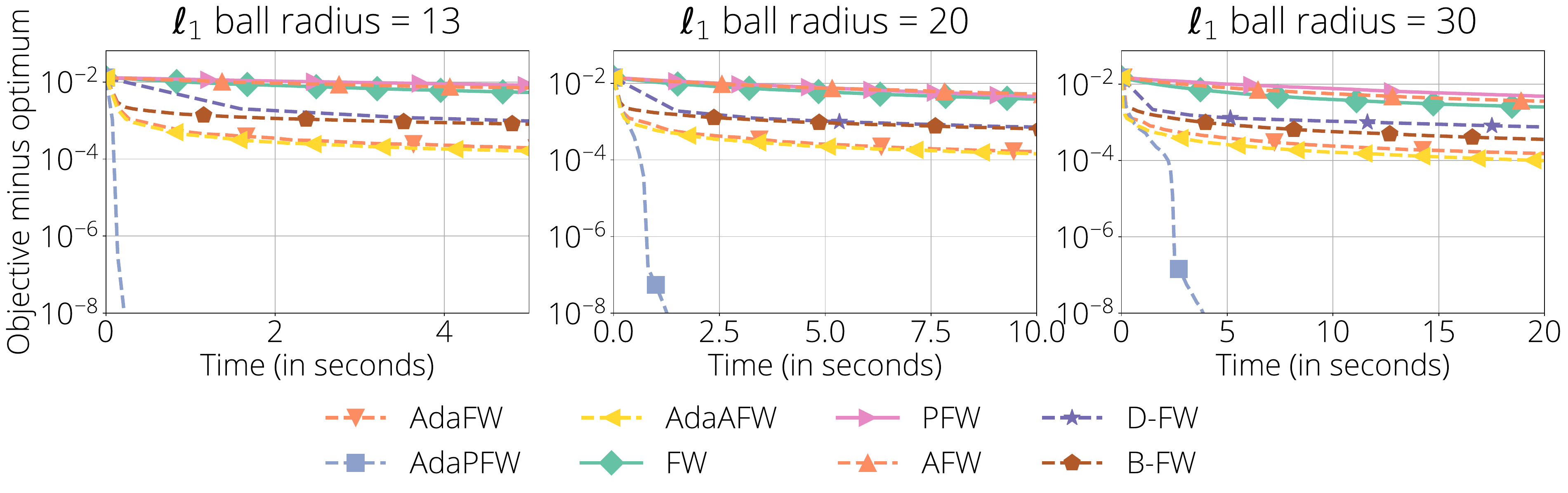}
  \caption{{\bfseries Comparison of different FW variants}. Problem is $\ell_1$-regularized logistic regression and dataset is Madelon in the first, RCV1 in the second figure. }
\end{figure}

\subsection{$\ell_1$-regularized logistic regression, RCV1 dataset}

The second experiment is identical to the first one, except the madelon datset is replaced by the larger RCV1 datset. Below we display the results of the comparison in this dataset:

\begin{figure}[ht]
  \centering\includegraphics[width=\linewidth]{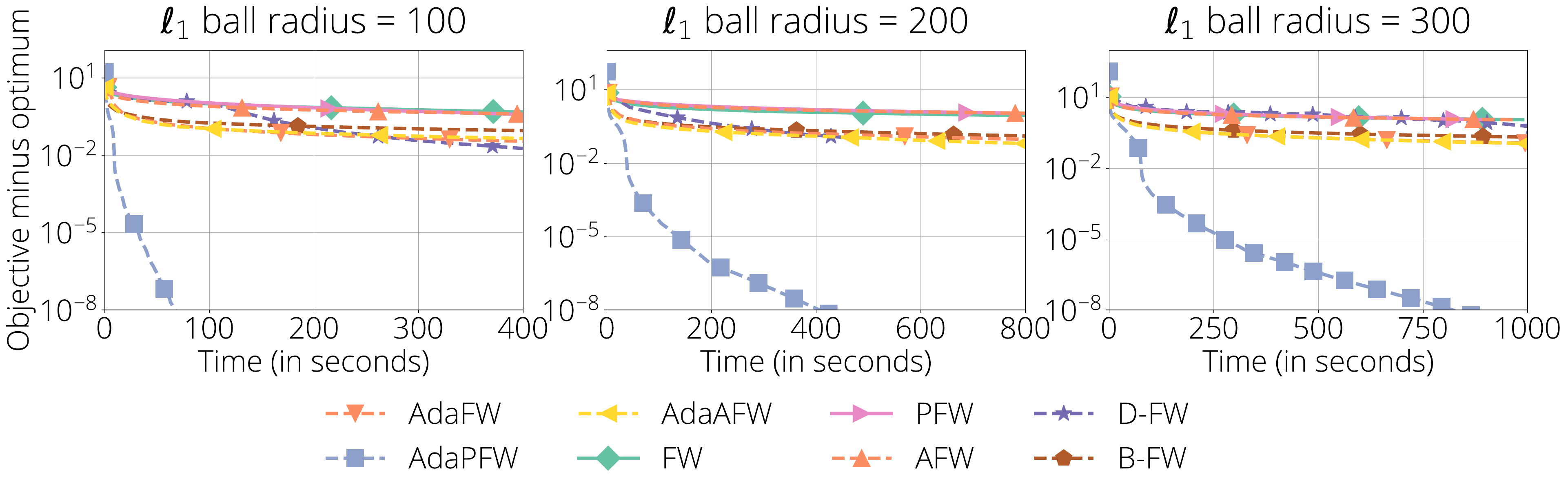}
  \caption{{\bfseries Comparison of different FW variants}. Problem is $\ell_1$-regularized logistic regression and dataset is RCV1. }
\end{figure}

\clearpage

\subsection{Nuclear norm-regularized Huber regression, MovieLens dataset}

For the third experiment, we consider a collaborative filtering problem with the Movielens 1M dataset~\citep{harper2016movielens} as provided by the spotlight\footnote{\url{https://github.com/maciejkula/spotlight}} Python package.

In this case the dataset consists of a sparse matrix $\boldsymbol{A}$ representing the ratings for the different movies and users. We denote by  $\mathcal{I}$ the non-zero indices of this matrix. Then the optimization probllem that we consider is the following
\begin{equation}
  \argmin_{\|\XX\|_* \leq \beta} \frac{1}{n} \sum_{(i, j) \in \mathcal{I}}^n L_\xi(\boldsymbol{A}_{i, j} -  \XX_{i, j})~,
\end{equation}
where $L_\xi$ is the Huber loss, defined as
\begin{equation}
  L_\xi (a) = \begin{cases}
    \frac{1}{2}{a^2}            & \text{for } |a| \le \xi, \\
    \xi (|a| - \frac{1}{2}\xi), & \text{otherwise\,.}
  \end{cases}
\end{equation}
The Huber loss is a quadratic for $|a|\leq \xi$ and grows linearly for $|a|>\xi$. The parameter $\xi$ controls this tradeoff and was set to 1 during the experiments.

We compared the variant of FW that do not require to store the active set on this problem (as these are the only competitive variants for this problem).
\begin{figure}[ht]
  \centering\includegraphics[width=0.7\linewidth]{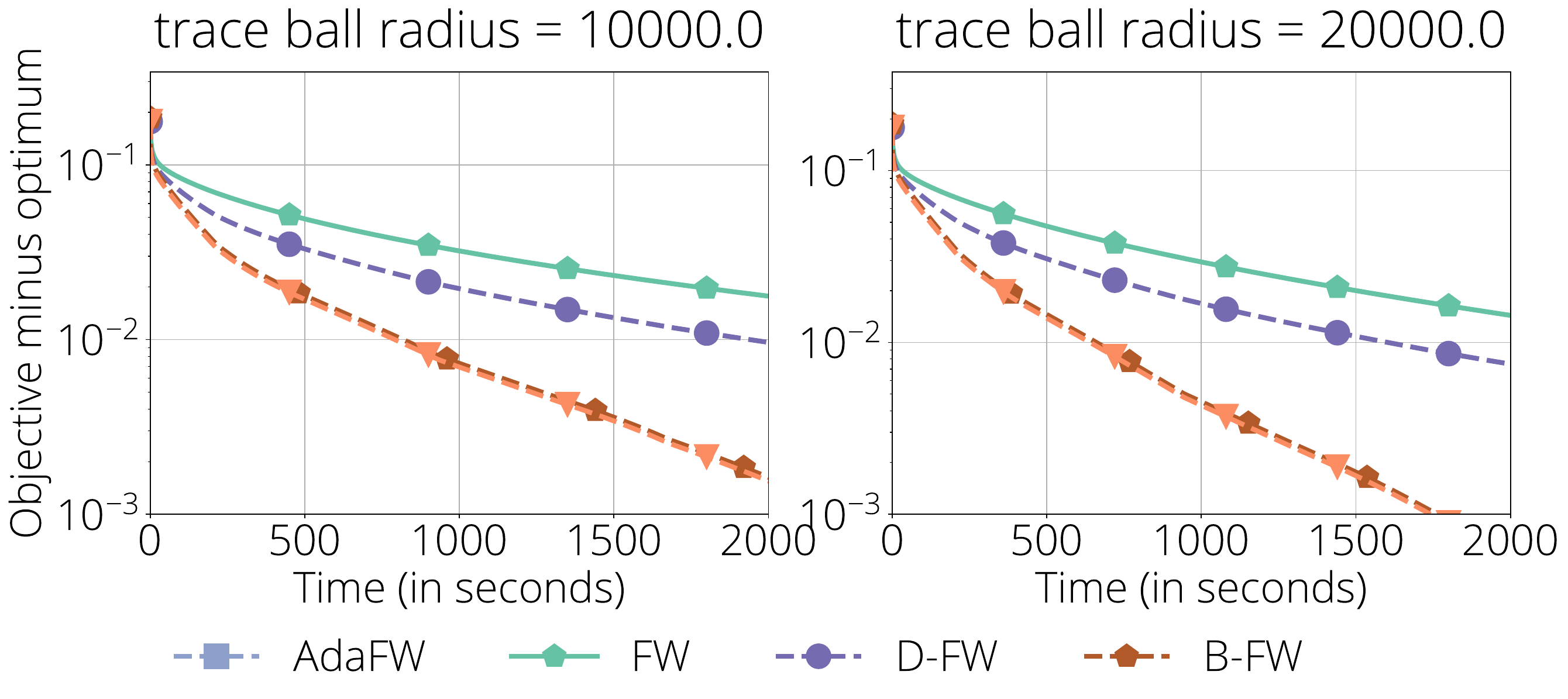}
  \caption{{\bfseries Comparison of different FW variants}. Comparison of FW variants on the Movielens 1M dataset. }
\end{figure}
\end{document}

%% file: defs.tex
\usepackage{nccmath} 
\definecolor{mydarkblue}{rgb}{0,0.08,0.45}
\definecolor{myblue}{HTML}{D2E4FC}
\newcommand*\mybluebox[1]{%
\colorbox{myblue}{\hspace{1em}#1\hspace{1em}}}

\definecolor{mygreen}{HTML}{008000}

\definecolor{myorange}{HTML}{FFC48C}

\definecolor{Gray}{gray}{0.92}

\newcommand{\red}{\color{red}}

\DeclareMathOperator*{\argmin}{\mathbf{arg\,min}}
\DeclareMathOperator*{\argmax}{\mathbf{arg\,max}}

\DeclareMathOperator*{\minimize}{minimize}

\DeclareMathOperator*{\diam}{diam}
\DeclareMathOperator*{\PWidth}{PWidth}
\DeclareMathOperator*{\mDW}{mDW}
\DeclareMathOperator*{\PdirW}{PdirW}
\DeclareMathOperator*{\radius}{radius}

\DeclareMathOperator*{\conv}{conv}
\DeclareMathOperator*{\face}{faces}
\DeclareMathOperator*{\cone}{cone}
\DeclareMathOperator*{\lin}{lin}

\usepackage{pifont}
\newcommand{\cmark}{\ding{51}}%
\newcommand{\xmark}{\ding{55}}%

\def\RR{{\mathbb R}}

\def\rr{{\boldsymbol r}}
\def\ss{{\boldsymbol s}}
\def\xx{{\boldsymbol x}}
\def\yy{{\boldsymbol y}}
\def\XX{{\boldsymbol X}}

\def\ss{{\boldsymbol s}}
\def\dd{{\boldsymbol d}}
\def\uu{{\boldsymbol u}}

\def\rr{{\boldsymbol r}}

\def\yy{{\boldsymbol y}}
\def\vv{{\boldsymbol v}}

\def\zz{{\boldsymbol z}}

\def\AA{{\mathcal A}}
\def\BB{{\mathcal B}}
\def\SS{{\mathcal S}}

\renewcommand{\aa}{\boldsymbol{a}}

\def\balpha{{\boldsymbol \alpha}}
\def\defas{\stackrel{\text{def}}{=}}

\newtheorem{theorem}{Theorem}
\newtheorem{definition}{Definition}
\newtheorem{lemma}{Lemma}
 \newtheorem{proposition}{Proposition}

\newtheorem{innercustomtheorem}{Theorem}
\newenvironment{customtheorem}[1]
  {\renewcommand\theinnercustomtheorem{#1}\innercustomtheorem}
  {\endinnercustomtheorem}


\newcounter{parentnumber}